\documentclass[12pt]{article}
\usepackage{latexsym}
\usepackage{amsmath}
\usepackage{amssymb}
\usepackage{amsfonts}
\usepackage{graphics,graphicx}
\usepackage{mathrsfs}
\usepackage{cite}
\usepackage{mathrsfs}
\usepackage{color}
\usepackage[colorlinks=true,linkcolor=blue]{hyperref}
\usepackage{tikz}

\numberwithin{equation}{section}

\newtheorem{theorem}{Theorem}[section]
\newtheorem{lemma}[theorem]{Lemma}
\newtheorem{pr}{Proposition}[section]

\newtheorem{definition}{Definition}[section]
\newtheorem{remark}[theorem]{Remark}

\newcommand{\eproof}{{\mbox{\ }~\hfill
\mbox{\large $\Box$} \par \vskip 10pt}}
\newcommand{\pf}{\noindent{\bf Proof}}
\newcommand{\eps}{\varepsilon}
\renewcommand{\l}{\ell}
\newcommand{\vph}{\varphi}
\newcommand{\C}{{\mathbb C}}
\newcommand{\re}{{\mbox{\rm Re\,}}}
\newcommand{\im}{{\mbox{{\rm Im\,}}}}

\newcommand{\R}{\mathbb R}

\title{Carleman estimate for complex second order elliptic operators with discontinuous Lipschitz coefficients}
\author{E.~Francini\thanks{Universit\`a di Firenze, Italy. Email: elisa.francini@unifi.it}
\qquad S.~Vessella\thanks{Universit\`a di Firenze, Italy. Email: sergio.vessella@unifi.it}
\qquad J-N.~Wang\thanks{National Taiwan University, Taiwan. Email: jnwang@ntu.edu.tw}}

\date{\today}

\begin{document}
\maketitle

\begin{abstract}
In this paper, we derive a local Carleman estimate for the complex second order elliptic operator with Lipschitz coefficients having jump discontinuities. Combing the result in \cite{BL} and the arguments in \cite{dflvw}, we present an elementary method to derive the Carleman estimate under the optimal regularity assumption on the coefficients. 
\end{abstract}

\section{Introduction}

Carleman estimates are important tools for proving the unique continuation property for partial differential equations. Additionally,  Carleman estimates have been successfully applied to study inverse problems and controllability of partial differential equations.  Most of Carleman estimates are proved under the assumption that the leading coefficients possess certain regularity. For example, for general second order elliptic operators, Carleman estimates were proved when the leading coefficients are at least Lipschitz \cite{ho3}. In general, the Lipschitz regularity assumption is the optimal condition for the unique continuation property to hold in $\R^n$ with $n\ge 3$ (see counterexamples constructed by Pli\'s \cite{pl} and Miller \cite{mi}). Therefore, Carleman estimates for second order elliptic operators with general discontinuous coefficients are most likely not valid. Nonetheless, recently, in the case of coefficients having jump discontinuities at an interface with homogeneous or non-homogeneous transmission conditions, one can still prove useful Carleman estimates, see, for example, Le Rousseau-Robbiano \cite{lr1}, \cite{lr2}, Le Rousseau-Lerner \cite{ll}, and \cite{dflvw}. 

Above mentioned results are proved for real coefficients. In many real world problems, the case of complex-valued coefficients arises naturally. The modeling of the current flows in biological tissues or the propagation of the electromagnetic waves in conductive media are typical examples. In these cases, the conductivities are complex-valued functions. On the other hand, in some situations, the conductivities are not continuous functions. For instance, in the human body, different organs have different conductivities. Therefore, to model the current flow in the human body, it is more reasonable to consider an anisotrotopic complex-valued conductivity with jump-type discontinuities \cite{mph06}.      

With potential applications in mind, our goal in this paper is to derive a Carleman estimate for the second order elliptic equations with complex-valued leading coefficients having jump-type discontinuities. Although such a Carleman estimate has been derived in \cite{BL}, we want to remark that the method used in \cite{BL}, also in \cite{lr1}, \cite{lr2}, and \cite{ll}, are based on the technique of pseudodifferential operators and hence requires $C^\infty$ coefficients and interface; while the method in \cite{dflvw} (and its parabolic counterpart, \cite{fv18}) relies on the Fourier transform and a version of partition of unity which requires only Lipschitz coefficients and $C^{1,1}$ interface. Hence, the main purpose of the paper is to extend the method in \cite{dflvw}, \cite{fv18} to second order elliptic operators with complex-valued coefficients. It is important to point out that even though second order elliptic operators with complex-valued coefficients can be written as a coupled second order elliptic system with real coefficients, neither the method in \cite{lr1}, \cite{lr2}, \cite{ll} nor that in \cite{dflvw} can be applied to coupled elliptic systems. Therefore, we need to work on operators with complex-valued coefficients directly. 

Our strategy to derive the Carleman estimate consists of two major steps. In the first step, we treat second order elliptic operators with constant complex coefficients. Based on \cite{BL}, by checking the strong pseudoconvexity and the transmission conditions in a neighborhood of a fixed point at the interface,  we can derive a Carleman estimate for second order elliptic operators with constant complex coefficients from \cite[Theorem~1.6]{BL}. We would like to mention that the result in \cite{BL} is stated for quite general complex coefficients, but here we can only verify the transmission condition with our choice of weight functions for complex coefficients having small imaginary parts. So in this paper we will consider this case. In the second step, we extend the Carleman estimate to the operator with non-constant complex coefficients with small imaginary parts. This method in this step is taken from the argument in \cite[Section 4]{dflvw}. The key tool is a version of partition of unity. 

Furthermore, in the second step, we need an interior Carleman estimate for second order elliptic operators having Lipschitz leading coefficients and with the weight function $\psi_\eps$. An interior Carleman estimate was proved in \cite[Theorem~8.3.1]{ho0}, but for operators with $C^1$ leading coefficients. Another interior estimate was established in \cite[Proposition~17.2.3]{ho3} for operators with Lipschitz leading coefficients, but with a different weight function. H\"ormander remarked in \cite{ho4} (page 703, line 7-8) that "Inspection of proof of Theorem~8.3.1 in \cite{ho0} shows that only Lipschitz continuity was actually used in the proof." But, as far as we can check, there is no formal proof of this statement in literature. To make the paper self contained, we would like give a detailed proof of interior Carleman estimate for second order elliptic operator with Lipschitz leading coefficients and with a rather general weight function, see Proposition~\ref{pr1}. This interior Carleman estimate may be useful on other occasions.  

In this paper, we present a detailed and elementary derivation of the Carleman estimate for the second order elliptic equations with complex-valued coefficients having jump-type discontinuities following our method in \cite{dflvw}. Having established the Carleman estimate, we then can apply the ideas in \cite{flvw} to prove a three-region inequality and those in \cite{CW} to prove a three-ball inequality across the interface. With the help of the three-ball inequality, we can study the size estimate problem for the complex conductivity equation following the ideas in \cite{CNW}. We will present these quantitative uniqueness results and the application to the size estimate in the forthcoming paper. 

The paper is organized as follows. In Section~\ref{pre}, we introduce notations that will be used in the paper and the statement of the theorem. In Section~\ref{sub-trans}, we derive a Carleman estimate for the operator having discontinuous piecewise constant coefficients. This Carleman estimate is a special case of \cite[Theorem 1.6]{BL}. Therefore, the main task of Section~\ref{sub-trans} is to check the transmission condition and the strong pseudoconvexity condition. Finally, the main Carleman estimate is proved in Section~\ref{sec4}. The key ingredient is a partition of unity introduced in \cite{dflvw}.

\section{Notations and statement of the main theorem}\label{pre}

We will state and prove the Carleman estimate for the case where the interface is flat. Since our Carleman estimate is local near any point at the interface, for general $C^{1,1}$ interface, it can be flatten by a suitable change of coordinates. Moreover, the transformed coefficients away from the interface remain Lipschitz. Define $H_{\pm}=\chi_{\mathbb{R}^n_{\pm}}$ where $\mathbb{R}^n_{\pm}=\{(x',x_n)\in \mathbb{R}^{n-1}\times\mathbb{R}|x_n\gtrless0\}$ and $\chi_{\mathbb{R}^n_{\pm}}$ is the characteristic function of $\mathbb{R}^n_{\pm}$.
In places we will use equivalently the symbols $\partial$, $\nabla$ and $D=-i\nabla$
to denote the gradient of a function and we will add the index $x'$ or $x_n$
to denote gradient in $\mathbb R^{n-1}$ and the derivative with respect to $x_n$
respectively. We further denote $\partial_\ell=\partial/\partial x_\ell$, $D_\ell=-i\partial_\ell$, and $\partial_{\xi_\ell}=\partial/\partial {\xi_\ell}$. 

Let $u_\pm\in C^\infty(\R^n)$. We define
\begin{equation*}\label{1.030}
u=H_+u_++H_-u_-=\sum_\pm H_{\pm}u_{\pm},
\end{equation*}
hereafter, we denote $\sum_\pm a_\pm=a_++a_-$,  and 
\begin{equation}\label{7.1}
\mathcal{L}(x,D)u:=\sum_{\pm}H_{\pm}{\rm div}(A_{\pm}(x)\nabla u_{\pm}),
\end{equation}
where
\begin{equation}\label{7.2}
A_{\pm}(x)=\{a^{\pm}_{\l j}(x)\}^n_{\l ,j=1}=\{a^{\pm}_{\l j}(x',x_n)\}^n_{\l ,j=1},\quad x'\in \mathbb{R}^{n-1},x_n\in \mathbb{R}
\end{equation}
is a Lipschitz symmetric matrix-valued function. Assume that
\begin{equation}\label{symm0}
a_{\l j}^{\pm}(x)=a_{j\l}^{\pm}(x),\quad\forall\;\;\l, j=1,\cdots,n,
\end{equation}
and furthermore
\begin{equation}\label{complex0}
a_{\l j}^{\pm}(x)=M_{\l j}^{\pm}(x)+i \gamma N_{\l j}^{\pm}(x),
\end{equation}
where $(M_{\l j}^{\pm})$ and $(N_{\l j}^{\pm})$ are real-valued matrices and $\gamma> 0$. We further assume that there exist $\lambda_0, \Lambda_0>0$ such that for all $\xi\in\R^n$ and $x\in \R^n$ we have
\begin{equation}\label{elliptic10}
\lambda_0|\xi|^2\le M^{\pm}(x)\xi\cdot\xi\le\Lambda_0|\xi|^2
\end{equation}
and
\begin{equation}\label{elliptic20}
\lambda_0|\xi|^2\le N^{\pm}(x)\xi\cdot\xi\le\Lambda_0|\xi|^2.
\end{equation}
In the paper, we consider Lipschitz coefficients $A_{\pm}$, i.e., there exists a constant $M_0>0$ such that
\begin{equation}\label{7.4}
|A_{\pm}(x)-A_{\pm}(y)|\leq M_0|x-y|.
\end{equation}
To treat the transmission conditions, we write
\begin{equation}\label{7.5}
h_0(x'):=u_+(x',0)-u_-(x',0),\ \forall\,\; x'\in \mathbb{R}^{n-1},
\end{equation}
\begin{equation}\label{7.6}
h_1(x'):=A_+(x',0)\nabla u_+(x',0)\cdot \nu-A_-(x',0)\nabla u_-(x',0)\cdot \nu,\ \forall\,\; x'\in \mathbb{R}^{n-1},
\end{equation}
where $\nu=e_n$.

Let us now introduce  the weight function. Let $\varphi$ be
\begin{equation}\label{2.1}
\varphi(x_n)=
\begin{cases}
\begin{array}{l}
\varphi_+(x_n):=\alpha_+x_n+\beta x_n^2/2,\quad x_n\geq 0,\\
\varphi_-(x_n):=\alpha_-x_n+\beta x_n^2/2,\quad x_n< 0,
\end{array}
\end{cases}
\end{equation}
where $\alpha_+$, $\alpha_-$ and $\beta$ are positive numbers which will be determined later. In what follows we denote by $\varphi_{+}$ and $\varphi_{-}$ the restriction of the weight function $\varphi$ to $[0,+\infty)$ and to $(-\infty,0)$ respectively. We use similar notation for any other weight functions. For any $\varepsilon>0$ let 
\begin{equation}\label{psi}
\psi_{\varepsilon}(x):=\varphi(x_n)-\frac{\varepsilon}{2}|x'|^2,
\end{equation}
and let 
\begin{equation}\label{wei}
\phi_{\delta}(x):=\psi_{\delta}(\delta^{-1}x),\quad\delta>0.
\end{equation}

For a function $h\in L^2(\mathbb{R}^{n})$, we define
\begin{equation*}
\hat{h}(\xi',x_n)=\int_{\mathbb{R}^{n-1}}h(x',x_n)e^{-ix'\cdot\xi}\,dx',\quad \xi'\in \mathbb{R}^{n-1}.
\end{equation*}
As usual we denote by $H^{1/2}(\mathbb{R}^{n-1})$ the space of the functions $f\in L^2(\mathbb{R}^{n-1})$ satisfying
$$\int_{\mathbb{R}^{n-1}}|\xi'||\hat{f}(\xi')|^2d\xi'<\infty,$$
with the norm
\begin{equation}\label{semR}
\|f\|^2_{H^{1/2}(\mathbb{R}^{n-1})}=\int_{\mathbb{R}^{n-1}}(1+|\xi'|^2)^{1/2}|\hat{f}(\xi')|^2d\xi'.
\end{equation}
Moreover we define
$$[f]_{1/2,\mathbb{R}^{n-1}}=\left[\int_{\mathbb{R}^{n-1}}\int_{\mathbb{R}^{n-1}}\frac{|f(x)-f(y)|^2}{|x-y|^n}dydx\right]^{1/2},$$
and recall that there is a positive constant $C$, depending only on $n$, such that
\begin{equation*}
C^{-1}\int_{\mathbb{R}^{n-1}}|\xi'||\hat{f}(\xi')|^2d\xi'\leq[f]^2_{1/2,\mathbb{R}^{n-1}}\leq C\int_{\mathbb{R}^{n-1}}|\xi'||\hat{f}(\xi')|^2d\xi',
\end{equation*}
so that the norm \eqref{semR} is equivalent to the norm $\|f\|_{L^2(\mathbb{R}^{n-1})}+[f]_{1/2,\mathbb{R}^{n-1}}$. We use the letters $C, C_0, C_1, \cdots$ to denote constants. The value of the constants may change from line to line, but it is always greater than $1$.

We will denote by $B'_r(x')$ the $(n-1)$-ball centered at $x'\in \mathbb{R}^{n-1}$ with radius $r>0$. Whenever $x'=0$ we denote $B'_r=B'_r(0)$. Likewise, we denote $B_r(x)$ be the $n$-ball centered at $x\in\R^n$ with radius $r>0$ and $B_r=B_r(0)$.
\bigskip

\begin{theorem}\label{thm8.2}
Let $u$ and $A_{\pm}(x)$ satisfy \eqref{7.1}-\eqref{7.6}. There exist $\alpha_+,\alpha_-,\beta, \delta_0, r_0, \gamma_0$ and $C$ depending on $\lambda_0, \Lambda_0, M_0$ such that if $\gamma<\gamma_0$, $\delta\le\delta_0$ and $\tau\geq C$, then
\begin{equation}\label{8.24}
\begin{aligned}
&\sum_{\pm}\sum_{k=0}^2\tau^{3-2k}\int_{\mathbb{R}^n_{\pm}}|D^k{u}_{\pm}|^2e^{2\tau\phi_{\delta,\pm}(x',x_n)}dx'dx_n+\sum_{\pm}\sum_{k=0}^1\tau^{3-2k}\int_{\mathbb{R}^{n-1}}|D^k{u}_{\pm}(x',0)|^2e^{2\phi_\delta(x',0)}dx'\\
&+\sum_{\pm}\tau^2[e^{\tau\phi_{\delta}(\cdot,0)}u_{\pm}(\cdot,0)]^2_{1/2,\mathbb{R}^{n-1}}+\sum_{\pm}[D(e^{\tau\phi_{\delta,\pm}}u_{\pm})(\cdot,0)]^2_{1/2,\mathbb{R}^{n-1}}\\
\leq &C\left(\sum_{\pm}\int_{\mathbb{R}^n_{\pm}}|\mathcal{L}(x,D)(u_{\pm})|^2\,e^{2\tau\phi_{\delta,\pm}(x',x_n)}dx'dx_n+[e^{\tau\phi_\delta(\cdot,0)}h_1]^2_{1/2,\mathbb{R}^{n-1}}\right.\\
&\left.+[D_{x'}(e^{\tau\phi_\delta}h_0)(\cdot,0)]^2_{1/2,\mathbb{R}^{n-1}}+\tau^{3}\int_{\mathbb{R}^{n-1}}|h_0|^2e^{2\tau\phi_\delta(x',0)}dx'+\tau\int_{\mathbb{R}^{n-1}}|h_1|^2e^{2\tau\phi_\delta(x',0)}dx'\right).
\end{aligned}
\end{equation}
where $u=H_+u_++H_-u_-$,  $u_{\pm}\in C^\infty(\mathbb{R}^{n})$ and ${\rm supp}\, u\subset B'_{\delta r_0}\times[-\delta r_0,\delta r_0]$, and $\phi_\delta$ is given by \eqref{wei}.
\end{theorem}
\begin{remark}
Estimate \eqref{8.24} is a local Carleman estimate near $x_n=0$. As mentioned above, by flattening the interface, we can derive a local Carleman estimate near a $C^{1,1}$ interface from \eqref{8.24}.  Nonetheless, an estimate like \eqref{8.24} is sufficient for some applications such as the inverse problem of estimating the size of an inclusion by one pair of boundary measurement {\rm (}see, for example, {\rm\cite{flvw}}{\rm )}. 
\end{remark}

\section{Carleman estimate for operators with constant coefficients}\label{sub-trans}

The purpose of this section is to derive \eqref{8.24} for ${\mathcal L}(x,D)$ with discontinuous piecewise constant coefficients. More precisely, we derive \eqref{8.24} for ${\mathcal L}_0(D)$, where ${\mathcal L}_0(D)$ is obtained from ${\mathcal L}(x,D)$ by freezing the variable $x$ at $(x'_0,0)$. Without loss of generality, we take $(x'_0,0)=(0,0)=0$ and thus
\[
{\mathcal L}_0(D)u=\mathcal{L}(0,D)u=\sum_{\pm}H_{\pm}{\rm div}(A_{\pm}(0)\nabla u_{\pm}).
\]
Since ${\mathcal L}_0$ has piecewise constant coefficients, to prove \eqref{8.24}, we will apply \cite[Theorem~1.6]{BL}. So the task here is to  verify the strong pseudoconvexity and transmission conditions for operator ${\mathcal L}_0$ with the weight function given in \eqref{psi}. 

To streamline the presentation,  we define $\Omega_1:=\{x_n<0\}, \Omega_2:=\{x_n>0\}$. On each side of the interface, we have complex second order elliptic operators. We denote
\[
P_k=\sum_{1\le j,\l\le n}a_{\l j}^{(k)}D_\ell D_j,\quad k=1,2,
\]
where $a_{\l j}^{(1)}=a_{\l j}^-$ and $a_{\l j}^{(2)}=a_{\l j}^+$. Here we denote $a_{\l j}^{(k)}=a_{\l j}^{(k)}(0)$. Corresponding to \eqref{symm0}-\eqref{elliptic20}, we have
\begin{equation}\label{symm}
a_{\l j}^{(k)}=a_{j\l}^{(k)}, 
\end{equation}
\begin{equation}\label{complex}
a_{\l j}^{(k)}=M_{\l j}^{(k)}+i \gamma N_{\l j}^{(k)},
\end{equation}
\begin{equation}\label{elliptic1}
\lambda_0|\xi|^2\le M^{(k)}\xi\cdot\xi\le\Lambda_0|\xi|^2,
\end{equation}
\begin{equation}\label{elliptic2}
\lambda_0|\xi|^2\le N^{(k)}\xi\cdot\xi\le\Lambda_0|\xi|^2.
\end{equation}
Since some computations in the verification of the transmission conditions are useful in proving the strong pseudoconvexity condition, we will begin with the discussion of the transmission conditions at the interface $\{x_n=0\}$.

\subsection{Transmission conditions}

We consider the natural transmission conditions that use the interface operators 
\[
T_k^1=(-1)^k,\quad T_k^2=(-1)^k\sum_{1\le j\le n}a_{nj}^{(k)}D_j
\]
that correspond to the continuity of the solution and of the normal flux, respectively. We now write the weight function
\begin{equation}\label{psieps}
\psi_\eps(x)=\vph(x_n)-\frac{\eps}{2}|x'|^2,
\end{equation}
where
\[
\vph(x_n)=\left\{
\begin{aligned}
&\vph_1(x_n),\quad x_n< 0\\
&\vph_2(x_n),\quad x_n\ge 0,
\end{aligned}\right.
\]
and
\[
\vph_k(x_n)=\alpha_kx_n+\frac 12\beta x_n^2
\]
with $\alpha_1, \alpha_2>0$ (corresponding to $\alpha_-$ and $\alpha_+$ in \eqref{2.1}, respectively) and $\beta>0$. Notice that $\vph$ is smooth in $\Omega_1, \Omega_2$ and is continuous across the interface.  Then we have
\[
\nabla\psi_\eps(0)=\left\{
\begin{aligned}
&(0,\cdots,0,\alpha_1),\quad x_n< 0\\
&(0,\cdots,0,\alpha_2),\quad x_n\ge 0.
\end{aligned}\right.
\]

Following the notations and the calculations in \cite[Section 1.7.1]{BL}, we have for $\omega:=(0,\xi',\nu,\tau)$ with $\xi'=(\xi_1,\cdots,\xi_{n-1})\ne 0, \nu=e_n$ and $\lambda\in\C$,
\[
{\tilde t}_{k,\psi_\eps}^1(\omega,\lambda)=(-1)^k
\]
and
\[
\begin{aligned}
{\tilde t}_{k,\psi_\eps}^2(\omega,\lambda)&=(-1)^ka_{nn}^{(k)}((-1)^k\lambda+i\tau\partial_{x_n}\psi_\eps(0))+(-1)^k\sum_{1\le j\le n-1}a_{nj}^{(k)}(\xi_j+i\tau\partial_{x_j}\psi_\eps(0))\\
&=(-1)^ka_{nn}^{(k)}((-1)^k\lambda+i\tau\alpha_k)+(-1)^k\sum_{1\le j\le n-1}a_{nj}^{(k)}\xi_j.
\end{aligned}
\]
The principal symbols of $P_k$, $k=1,2$, can be written as
\begin{equation}\label{pk}
p_k(\xi)=a_{nn}^{(k)}((\xi_n+\sum_{1\le j\le n-1}\frac{a_{nj}^{(k)}}{a_{nn}^{(k)}}\xi_j)^2+b_k(\xi')),
\end{equation}
where
\begin{equation}\label{bk}
b_k(\xi')=(a_{nn}^{(k)})^{-2}\sum_{1\le\l, j\le n-1}(a_{\l j}^{(k)}a_{nn}^{(k)}-a_{n\l}^{(k)}a_{nj}^{(k)})\xi_\l\xi_j.
\end{equation}
We also need to introduce the principal symbol of the conjugate operators
\begin{equation}\label{tpk}
\begin{aligned}
\tilde p_{k,\psi_\eps}(\omega,\lambda)&=a_{nn}^{(k)}\Big{[}\big{(}(-1)^k\lambda+i\tau\partial_{x_n}\psi_\eps(0)+\sum_{1\le j\le n-1}\frac{a_{nj}^{(k)}}{a_{nn}^{(k)}}(\xi_j+i\tau\partial_{x_j}\psi_\eps(0))\big{)}^2\\
&\qquad +b_k(\xi'+i\tau\partial_{x'}\psi_\eps(0))\Big{]}\\
&=a_{nn}^{(k)}\Big{[}\big{(}(-1)^k\lambda+i\tau\alpha_k+\sum_{1\le j\le n-1}\frac{a_{nj}^{(k)}}{a_{nn}^{(k)}}\xi_j\big{)}^2+b_k(\xi')\Big{]}.
\end{aligned}
\end{equation}

Let us introduce $A^{(k)}, B^{(k)}\in\R$ for $k=1,2$ such that
\begin{equation}\label{bk0}
\begin{aligned}
b_k(\xi')&=(a_{nn}^{(k)})^{-2}\sum_{1\le\l, j\le n-1}(a_{\l j}^{(k)}a_{nn}^{(k)}-a_{n\l}^{(k)}a_{nj}^{(k)})\xi_\l\xi_j\\
&=(A^{(k)}-iB^{(k)})^2,
\end{aligned}
\end{equation}
where $A^{(k)}\ge 0$. We also denote
\begin{equation}\label{efk}
\sum_{1\le j\le n-1}\frac{a_{nj}^{(k)}}{a_{nn}^{(k)}}\xi_j=E^{(k)}+iF^{(k)},
\end{equation}
where $E^{(k)}, F^{(k)}\in\R$. Using \eqref{tpk}, \eqref{bk0}, and \eqref{efk}, we can write
\[
\begin{aligned}
\tilde p_{2,\psi_\eps}&=a_{nn}^{(2)}[(\lambda+i\tau\alpha_2+E^{(2)}+iF^{(2)})^2+(A^{(2)}-iB^{(2)})^2]\\
&=a_{nn}^{(2)}[(\lambda+i\tau\alpha_2+E^{(2)}+iF^{(2)}+i(A^{(2)}-iB^{(2)}))\\
&\quad\cdot(\lambda+i\tau\alpha_2+E^{(2)}+iF^{(2)}-i(A^{(2)}-iB^{(2)}))]\\
&=a_{nn}^{(2)}(\lambda-\sigma_1^{(2)})(\lambda-\sigma_2^{(2)}),
\end{aligned}
\]
where
\[
\begin{aligned}
\sigma_1^{(2)}&=-E^{(2)}-B^{(2)}-i(\tau\alpha_2+F^{(2)}+A^{(2)}),\\
\sigma_2^{(2)}&=-E^{(2)}+B^{(2)}-i(\tau\alpha_2+F^{(2)}-A^{(2)}).
\end{aligned}
\]
On the other hand, we can write
\[
\begin{aligned}
\tilde p_{1,\psi_\eps}&=a_{nn}^{(1)}[(-\lambda+i\tau\alpha_1+E^{(1)}+iF^{(1)})^2+(A^{(1)}-iB^{(1)})^2]\\
&=a_{nn}^{(1)}[(\lambda-i\tau\alpha_1-E^{(1)}-iF^{(1)}+i(A^{(1)}-iB^{(1)}))\\
&\quad\cdot(\lambda-i\tau\alpha_1-E^{(1)}-iF^{(1)}-i(A^{(1)}-iB^{(1)}))]\\
&=a_{nn}^{(1)}(\lambda-\sigma_1^{(1)})(\lambda-\sigma_2^{(1)}),
\end{aligned}
\]
where
\[
\begin{aligned}
\sigma_1^{(1)}&=E^{(1)}+B^{(1)}+i(\tau\alpha_1+F^{(1)}+A^{(1)}),\\
\sigma_2^{(1)}&=E^{(1)}-B^{(1)}+i(\tau\alpha_1+F^{(1)}-A^{(1)}).
\end{aligned}
\]
Let us introduce the polynomial
\[
K_{k,\psi_\eps}(\omega,\lambda):=\prod_{{\scriptsize{\im}}\sigma_j^{(k)}\ge 0}(\lambda-\sigma_j^{(k)}).
\]
Now we state the definition of transmission conditions given in \cite[Definition~1.4]{BL}.
\begin{definition}\label{def3.1}
The pair $\{P_k,\psi_\eps, T_k^j,\; k=1,2,\; j=1,2\}$ satisfies the transmission condition at $\omega$ if for any polynomials $q_1(\lambda), q_2(\lambda)$, there exist polynomials $U_1(\lambda), U_2(\lambda)$ and constant $c_1, c_2$ such that
\[
\left\{
\begin{aligned}
q_1(\lambda)&=c_1\tilde t^1_{1,\psi_\eps}(\omega,\lambda)+c_2\tilde t^2_{1,\psi_\eps}(\omega,\lambda)+U_1(\lambda)K_{1,\psi_\eps}(\omega,\lambda),\\
q_2(\lambda)&=c_1\tilde t^1_{2,\psi_\eps}(\omega,\lambda)+c_2\tilde t^2_{2,\psi_\eps}(\omega,\lambda)+U_2(\lambda)K_{2,\psi_\eps}(\omega,\lambda).
\end{aligned}\right.
\]
\end{definition}

In order to check the transmission conditions, we need to study the polynomial $K_{k,\psi_\eps}(\omega,\lambda)$. For this reason, we need to determine the signs of the imaginary parts of the roots $\sigma_j^{(k)}$ defined above. Note that we can write
\begin{equation}\label{0517}
b_k(\xi')=\frac{1}{a_{nn}^{(k)}}\sum_{1\le\l, j\le n-1}a_{\l j}^{(k)}\xi_\l \xi_j-(E^{(k)}+iF^{(k)})^2.
\end{equation}
Since $b_k$ plays an essential role, we begin by working some calculations on the matrix $\frac{1}{a_{nn}^{(k)}}{\cal A}^{(k)}$, where ${\cal A}^{(k)}$ is the matrix $(a_{\l j}^{(k)})$. Let $a_{nn}^{(k)}=|a_{nn}^{(k)}|e^{i\theta}$. Choosing $\xi=e_n$, we have that
\[
a_{nn}^{(k)}=\sum_{1\le\l, j\le n}a_{\l j}^{(k)}\xi_\l\xi_j=\sum_{1\le\l, j\le n}M_{\l j}^{(k)}\xi_\l\xi_j+i\gamma\sum_{1\le\l, j\le n}N_{\l j}^{(k)}\xi_\l\xi_j.
\]
Hence, from \eqref{elliptic1}, \eqref{elliptic2}, we have that
\[
\lambda_0\le\re(a_{nn}^{(k)})\le\Lambda_0\quad\mbox{and}\quad\lambda_0\le\frac{\im(a_{nn}^{(k)})}{\gamma}\le\Lambda_0
\]
and so that $\theta\in[0,\pi/2)$. Let us evaluate
\begin{equation}\label{04-8}
\begin{aligned}
(a_{nn}^{(k)})^{-1}{\cal A}^{(k)}&=|a_{nn}^{(k)}|^{-1}(M^{(k)}+i\gamma N^{(k)})(\cos\theta-i\sin\theta)\\
&=|a_{nn}^{(k)}|^{-1}[\cos\theta M^{(k)}+\gamma \sin\theta N^{(k)}+i(-\sin\theta M^{(k)}+\gamma\cos\theta N^{(k)})].
\end{aligned}
\end{equation}
Using \eqref{elliptic1}, \eqref{elliptic2} again, we see that for $\xi\in\R^n$
\begin{equation}\label{55}
\begin{aligned}
\re((a_{nn}^{(k)})^{-1}{\cal A}^{(k)}\xi\cdot\xi)&=|a_{nn}^{(k)}|^{-1}[\cos\theta M^{(k)}\xi\cdot\xi+\gamma \sin\theta N^{(k)}\xi\cdot\xi]\\
&\ge|a_{nn}^{(k)}|^{-1}\lambda_0(\cos\theta+\gamma\sin\theta)|\xi|^2. 
\end{aligned}
\end{equation}
In fact, since $\cos\theta=M_{nn}^{(k)}|a_{nn}^{(k)}|^{-1}$ and $\sin\theta=\gamma N_{nn}^{(k)}|a_{nn}^{(k)}|^{-1}$, while $|a_{nn}^{(k)}|^2=(M_{nn}^{(k)})^2+\gamma^2(N_{nn}^{(k)})^2$, we have
\begin{equation}\label{66}
|a_{nn}^{(k)}|^{-1}(\cos\theta+\gamma\sin\theta)=\frac{M_{nn}^{(k)}+\gamma^2N_{nn}^{(k)}}{(M_{nn}^{(k)})^2+\gamma^2(N_{nn}^{(k)})^2}\ge\frac{\lambda_0(1+\gamma^2)}{\Lambda_0^2(1+\gamma^2)}=\frac{\lambda_0}{\Lambda_0^2}.
\end{equation}
Combining \eqref{55} and \eqref{66} implies
\begin{equation}\label{88}
\re((a_{nn}^{(k)})^{-1}{\cal A}^{(k)}\xi\cdot\xi)\ge\frac{\lambda_0^2}{\Lambda_0^2}|\xi|^2:=\tilde\lambda_1|\xi|^2.
\end{equation}

Now let us write
\begin{equation}\label{20}
\begin{aligned}
&\tilde\lambda_1|\xi|^2\le\re((a_{nn}^{(k)})^{-1}{\cal A}^{(k)}\xi\cdot\xi)\\
=&\re[\sum_{1\le\l, j\le n-1}\frac{a_{\l j}^{(k)}}{a_{nn}^{(k)}}\xi_\l\xi_j+2\sum_{1\le j\le n-1}\frac{a_{nj}^{(k)}}{a_{nn}^{(k)}}\xi_n\xi_j+\xi_n^2]\\
=&\xi_n^2+2b_0^{(k)}(\xi')\xi_n+b_1^{(k)}(\xi'),
\end{aligned}
\end{equation}
where
\[
b_0^{(k)}(\xi')=\re(\sum_{1\le j\le n-1}\frac{a_{nj}^{(k)}}{a_{nn}^{(k)}}\xi_j)=\re(E^{(k)}+iF^{(k)})=E^{(k)}
\]
and
\[
b_1^{(k)}(\xi')=\re(\sum_{1\le\l, j\le n-1}\frac{a_{\l j}^{(k)}}{a_{nn}^{(k)}}\xi_\l\xi_j).
\]
Substituting $\tilde\xi_n=\xi_n=-b_0^{(k)}(\xi')$ into \eqref{20} gives
\[
\tilde\lambda_1(|\xi'|^2+|\tilde\xi_n|^2)\le\tilde\xi_n^2-2b_0^{(k)}(\xi')\tilde\xi_n+b_1^{(k)}(\xi')=-(b_0^{(k)}(\xi')^2+b_1^{(k)}(\xi'),
\]
which implies
\begin{equation}\label{50}
\tilde\lambda_1|\xi'|^2\le\re(\sum_{1\le\l, j\le n-1}\frac{a_{\l j}^{(k)}}{a_{nn}^{(k)}}\xi_\l\xi_j)-E_k^2.
\end{equation}
Putting \eqref{0517} and \eqref{50} together gives
\begin{equation}\label{501}
\begin{aligned}
\re(b_k(x_0,\xi'))&=\re(\sum_{1\le\l, j\le n-1}\frac{a_{\l j}^{(k)}}{a_{nn}^{(k)}}\xi_\l\xi_j)-(E^{(k)})^2+(F^{(k)})^2\\
&\ge\tilde\lambda_1|\xi'|^2+(F^{(k)})^2>0.
\end{aligned}
\end{equation}

The following lemma guarantees the positivity of $A^{(k)}$.
\begin{lemma}\label{ak}
Assume that \eqref{elliptic1} and \eqref{elliptic2} hold. Then 
\begin{equation}\label{akk}
A^{(k)}\ge\sqrt{\tilde\lambda_1|\xi'|^2+|F^{(k)}|^2}>|F^{(k)}|.
\end{equation}
\end{lemma}
\pf From \eqref{bk0}, it is easy to see that
\[
A^{(k)}=\re\sqrt{b_k}=\sqrt{\frac{a+\sqrt{a^2+b^2}}{2}},
\]
where $a=\re b_k$ and $b=\im b_k$. We have from \eqref{501} that $a>0$ and thus
\[
A^{(k)}\ge\sqrt{a}\ge\sqrt{\tilde\lambda_1|\xi'|^2+(F^{(k)})^2}>|F^{(k)}|.
\]
\eproof

Lemma~\ref{ak} implies
\begin{equation}\label{001}
\begin{aligned}
\im\sigma_1^{(2)}&=-(\tau\alpha_2+F^{(2)}+A^{(2)})=-\tau\alpha_2-F^{(2)}-A^{(2)}\\
&\le-\tau\alpha_2-|F^{(2)}|-F^{(2)}\le-\tau\alpha_2<0
\end{aligned}
\end{equation}
and
\begin{equation}\label{002}
\im\sigma_1^{(1)}=\tau\alpha_1+F^{(1)}+A^{(1)}>\tau\alpha_1+F^{(1)}+|F^{(1)}|\ge\tau\alpha_1>0.
\end{equation}

We are now ready to check the transmission condition defined in Definition~\ref{def3.1}. Being able to satisfy this condition depends on the degree of $K_{1,\psi_\eps}$ and $K_{2,\psi_\eps}$, that is, on the number of roots with negative imaginary parts. 

\bigskip\noindent
{\bf Case 1}. $\tilde p_{2,\psi_\eps}$ has two roots in $\{\im z<0\}$, i.e., $-\tau\alpha_2-F^{(2)}+A^{(2)}<0$ in view of \eqref{001}. In this case, we have that
\[
K_{2,\psi_\eps}=1,\;\;\mbox{\rm while}\;\; K_{1,\psi_\eps}\;\;\mbox{\rm has degree}\; 1\;\mbox{or}\; 2\;\; (\mbox{note}\; \eqref{002}). 
\]
Since ${\tilde t}_{2,\psi_\eps}^1(\omega,\lambda)=1$ and
\[
{\tilde t}_{2,\psi_\eps}^2(\omega,\lambda)=a_{nn}^{(2)}(\lambda+i\tau\alpha_2+\sum_{1\le j\le n-1}\frac{a_{nj}^{(2)}}{a_{nn}^{(2)}}\xi_j),
\]
for any $q_2(\lambda)$, we simply choose
\[
U_2(\lambda)=q_2(\lambda)-c_1\tilde t^1_{2,\psi_\eps}-c_2\tilde t^2_{2,\psi_\eps}.
\]
On the other hand, we have ${\tilde t}_{1,\psi_\eps}^1(\omega,\lambda)=-1$ and
\[
{\tilde t}_{1,\psi_\eps}^2(\omega,\lambda)=a_{nn}^{(1)}(\lambda-i\tau\alpha_1-\sum_{1\le j\le n-1}\frac{a_{nj}^{(1)}}{a_{nn}^{(1)}}\xi_j).
\]
Then for any polynomial $q_1(\lambda)$, we choose $U_1(\lambda)$ to be the quotient of the division between $q_1$ and $K_{1,\psi_\eps}$. The remainder term is equal to $c_1\tilde t_{1,\psi_\eps}+c_2\tilde t_{2,\psi_\eps}$ with suitable $c_1$, $c_2$.

\bigskip\noindent
{\bf Case 2}. Assume that $\im \sigma_2^{(2)}\ge 0$ and $\im \sigma_2^{(1)}\ge 0$, i.e.,
\[
-\tau\alpha_2-F^{(2)}+A^{(2)}\ge 0,\quad\tau\alpha_1+F^{(1)}-A^{(1)}\ge 0.
\]
Then $K_{1,\psi_\eps}$ has degree $2$ and $K_{2,\psi_\eps}$ has degree $1$. In order to avoid this case, we need to be sure that if $-\tau\alpha_2-F^{(2)}+A^{(2)}\ge 0$, then $\tau\alpha_1+F^{(1)}-A^{(1)}<0$, that is,
\[
\tau\alpha_2+F^{(2)}-A^{(2)}\le 0\Rightarrow \tau\alpha_1+F^{(1)}-A^{(1)}<0.
\]
This can be achieved by assuming that
\begin{equation}\label{tc20}
\frac{\alpha_2}{\alpha_1}>\frac{A^{(2)}-F^{(2)}}{A^{(1)}-F^{(1)}},\;\;\forall\;\;\xi'\ne 0.
\end{equation}
Recall that $A^{(k)}-F^{(k)}>0$, $k=1,2$. We remark that all $A^{(k)}$ and $F^{(k)}$ are homogeneous of degree $1$ in $\xi'$.  Hence \eqref{tc20} holds provided
\begin{equation}\label{tc2}
\frac{\alpha_2}{\alpha_1}=\max_{|\xi'|=1}\left\{\frac{A^{(2)}-F^{(2)}}{A^{(1)}-F^{(1)}}\right\}+1.
\end{equation}
Hence, if we assume \eqref{tc2}, then the transmission condition is satisfied. 

\bigskip\noindent
{\bf Case 3}. Each symbol has exactly one root in $\{\im z<0\}$, i.e.,
\[
\tau\alpha_1+F^{(1)}-A^{(1)}<0,\quad -\tau\alpha_2-F^{(2)}+A^{(2)}>0.
\]
In this case, we have
\[
K_{1,\psi_\eps}=(\lambda-\sigma_1^{(1)}),\quad K_{2,\psi_\eps}=(\lambda-\sigma_2^{(2)}).
\]
Given polynomials $q_1(\lambda), q_2(\lambda)$, there exist $U_1(\lambda), U_2(\lambda)$ such that 
\[
\begin{aligned}
q_1(\lambda)&=U_1(\lambda)K_{1,\psi_\eps}+\tilde q_1,\\
q_2(\lambda)&=U_2(\lambda)K_{2,\psi_\eps}+\tilde q_2,
\end{aligned}
\]
where $\tilde q_1, \tilde q_2$ are constants in $\lambda$. The transmission condition is satisfied if there exists constants $\mu_1, \mu_2, c_1, c_2$ so that
\[
\left\{\begin{aligned}
\tilde q_1&=\mu_1K_{1,\psi_\eps}+c_1\tilde t_{1,\psi_\eps}^1+c_2\tilde t_{1,\psi_\eps}^2,\\
\tilde q_2&=\mu_2K_{2,\psi_\eps}+c_1\tilde t_{2,\psi_\eps}^1+c_2\tilde t_{2,\psi_\eps}^2,
\end{aligned}\right.
\]
namely,
\begin{equation}\label{tc0}
\left\{\begin{aligned}
\tilde q_1&=\mu_1(\lambda-\sigma_1^{(1)})-c_1+c_2a_{nn}^{(1)}(\lambda-i\tau\alpha_1-E^{(1)}-iF^{(1)})\\
\tilde q_2&=\mu_2(\lambda-\sigma_2^{(2)})+c_1+c_2a_{nn}^{(2)}(\lambda+i\tau\alpha_2+E^{(2)}+iF^{(2)}).
\end{aligned}\right.
\end{equation}
System \eqref{tc0} is equivalent to 
\begin{equation}\label{tc1}
\left\{
\begin{aligned}
&\mu_1+c_2a_{nn}^{(1)}=0\\
&\mu_2+c_2a_{nn}^{(2)}=0\\
&\mu_1\sigma_1^{(1)}+c_1+c_2a_{nn}^{(1)}(i\tau\alpha_1+E^{(1)}+iF^{(1)})=-\tilde q_1\\
&-\mu_2\sigma_2^{(2)}+c_1+c_2a_{nn}^{(2)}(i\tau\alpha_2+E^{(2)}+iF^{(2)})=\tilde q_2.
\end{aligned}\right.
\end{equation}

System \eqref{tc1} has a unique solution if and only if the matrix
\[
T=\begin{pmatrix}1&0&0&a_{nn}^{(1)}\\0&1&0&a_{nn}^{(2)}\\\sigma_1^{(1)}&0&1&\zeta_1\\0&-\sigma_2^{(2)}&1&\zeta_2\end{pmatrix}
\]
with $\zeta_1=a_{nn}^{(1)}(i\tau\alpha_1+E^{(1)}+iF^{(1)})$, $\zeta_2=a_{nn}^{(2)}(i\tau\alpha_2+E^{(2)}+iF^{(2)})$, is nonsingular. We compute
\[
\begin{aligned}
\mbox{\rm det}T=&\mbox{\rm det}\begin{pmatrix}1&0&a_{nn}^{(1)}\\0&1&a_{nn}^{(2)}\\0&-\sigma_2^{(2)}&\zeta_2\end{pmatrix}-\mbox{\rm det}\begin{pmatrix}1&0&a_{nn}^{(1)}\\0&1&a_{nn}^{(2)}\\\sigma_1^{(1)}&0&\zeta_1\end{pmatrix}\\
=&\zeta_2+\sigma_2^{(2)}a_{nn}^{(2)}-\zeta_1+\sigma_1^{(1)}a_{nn}^{(1)}\\
=&a_{nn}^{(2)}(i\tau\alpha_2+E^{(2)}+iF^{(2)}-E^{(2)}+B^{(2)}-i\tau\alpha_2-iF^{(2)}+iA^{(2)})\\
&+a_{nn}^{(1)}(-i\tau\alpha_1-E^{(1)}-iF^{(1)}+E^{(1)}+B^{(1)}+i\tau\alpha_1+iF^{(1)}+iA^{(1)})\\
=&a_{nn}^{(2)}(B^{(2)}+iA^{(2)})+a_{nn}^{(1)}(B^{(1)}+iA^{(1)}).
\end{aligned}
\]
Therefore, if 
\begin{equation}\label{ne0}
a_{nn}^{(2)}(B^{(2)}+iA^{(2)})+a_{nn}^{(1)}(B^{(1)}+iA^{(1)})\ne 0,
\end{equation}
then the transmission condition holds.

We now verify \eqref{ne0}. In the real case where $a_{nn}^{(2)}, a_{nn}^{(1)}$ are positive real numbers, it is easy to see that
\[
a_{nn}^{(2)}A^{(2)}+a_{nn}^{(1)}A^{(1)}>0
\]
and thus \eqref{ne0} holds. 

For the complex case, we want to show that there exists $\gamma_0>0$ such that if $\gamma<\gamma_0$, then \eqref{ne0} is satisfied. Let $u_k=A^{(k)}+i B^{(k)}$ and $v_k=iu_k=-B^{(k)}+iA^{(k)}$. We will consider $u_k$ and $v_k$ as vectors in $\R^2$, i.e., $u_k=(A^{(k)},B^{(k)})$, $v_k=u_k^\perp=(-B^{(k)},A^{(k)})$. Let $a_{nn}^{(k)}=\eta^{(k)}+i\gamma\delta^{(k)}$ for $\eta^{(k)}, \delta^{(k)}\in\R$. By the ellipticity conditions \eqref{elliptic1}, \eqref{elliptic2}, we have
\[
\lambda_0\le\eta^{(k)}\le\Lambda_0,\quad\lambda_0\le\delta^{(k)}\le\Lambda_0.
\]
Notice that $\mbox{\rm det}T=0$ if and only if
\[
(\eta^{(2)}+i\gamma\delta^{(2)})(B^{(2)}+iA^{(2)})+(\eta^{(1)}+i\gamma\delta^{(1)})(B^{(1)}+iA^{(1)})=0,
\]
i.e.,
\[
\begin{aligned}
&(\eta^{(2)}B^{(2)}-\gamma\delta^{(2)}A^{(2)}+\eta^{(1)}B^{(1)}-\gamma\delta^{(1)}A^{(1)})\\
&\;\;+i(\eta^{(2)}A^{(2)}+\gamma\delta^{(2)}B^{(2)}+\eta^{(1)}A^{(1)}+\gamma\delta^{(1)}B^{(1)})=0,
\end{aligned}
\]
which is equivalent to
\begin{equation}\label{tc5}
\eta^{(2)}\begin{pmatrix}A^{(2)}\\B^{(2)}\end{pmatrix}+\eta^{(1)}\begin{pmatrix}A^{(1)}\\B^{(1)}\end{pmatrix}=\gamma\delta^{(2)}\begin{pmatrix}-B^{(2)}\\A^{(2)}\end{pmatrix}+\gamma\delta^{(1)}\begin{pmatrix}-B^{(1)}\\A^{(1)}\end{pmatrix}
\end{equation}
or simply
\begin{equation}\label{tc6}
\eta^{(2)}u_2+\eta^{(1)}u_1=\gamma\delta^{(2)}v_2+\gamma\delta^{(1)}v_1.
\end{equation}

Recall that $A^{(k)}\ge |F^{(k)}|>0$. Therefore, in the real case $\gamma\delta^{(k)}=0$, then \eqref{tc5} will never be satisfied. If $B^{(1)}$ and $B^{(2)}$ have the same sign, that is, either $B^{(k)}\ge 0$ or $B^{(k)}\le 0$ for $k=1,2$, \eqref{tc6} can not hold. To see this, let us consider $B^{(k)}\ge 0$, $k=1,2$. Then $u_1, u_2$ are in the first quadrant of the plane and $v_1, v_2$ are in the second quadrant of the plane. The sets
\[
C_u=\{\eta^{(2)}u_2+\eta^{(1)}u_1: \eta^{(k)}\ge 0\},\quad C_u=\{\gamma\delta^{(2)}v_2+\gamma\delta^{(1)}v_1: \gamma\delta^{(k)}\ge 0\}
\]
can only intersect at the original. Same thing happens if $B^{(k)}\le 0$ for $k=1,2$. 

The only case we need to investigate is when $B^{(1)}$ and $B^{(2)}$ have different signs. For example, let us assume
\[
B^{(1)}>0,\quad B^{(2)}<0.
\]
Even in this case, the intersection between $C_u$ and $C_v$ is non-trivial if  the angle $\phi$ between $u_1$ and $u_2$ is less than $\pi/2$. Note that $u_1$ is the first quadrant and $u_2$ is in the fourth quadrant. So the angle between $u_1$ and $u_2$ is less than $\pi$. We would like to show that \eqref{tc6} cannot hold for $\phi\in[\pi/2,\pi)$ if we choose $\gamma_0$ small enough. 

Note that in this case $\cos\phi\le 0$. To do so, we estimate $\|\eta^{(2)}u_2+\eta^{(1)}u_1\|$ from below and $\|\delta^{(2)}v_2+\delta^{(1)}v_1\|$ from above. We now discuss the estimate of $\|\delta^{(2)}v_2+\delta^{(1)}v_1\|$ from above. Compute
\begin{equation}\label{a1}
\begin{aligned}
\|\delta^{(2)}v_2+\delta^{(1)}v_1\|^2&=(\delta^{(2)})^2[(A^{(2)})^2+(B^{(2)})^2]+(\delta^{(1)})^2[(A^{(1)})^2+(B^{(1)})^2]\\
&+2\delta^{(1)}\delta^{(2)}(-B^{(2)},A^{(2)})\cdot(-B^{(1)},A^{(1)})\\
&=(\delta^{(2)})^2[(A^{(2)})^2+(B^{(2)})^2]+(\delta^{(1)})^2[(A^{(1)})^2+(B^{(1)})^2]\\
&+2\delta^{(1)}\delta^{(2)}[(A^{(2)})^2+(B^{(2)})^2]^{1/2}[(A^{(1)})^2+(B^{(1)})^2]^{1/2}\cos\phi\\
&\le (\delta^{(2)})^2[(A^{(2)})^2+(B^{(2)})^2]+(\delta^{(1)})^2[(A^{(1)})^2+(B^{(1)})^2].
\end{aligned}
\end{equation}
In view of \eqref{bk0} and \eqref{0517}, we have
\begin{equation}\label{a2}
\begin{aligned}
(A^{(k)})^2+(B^{(k)})^2=|b_k|&=|\sum_{1\le\l, j\le n-1}\frac{a_{\l j}^{(k)}}{a_{nn}^{(k)}}\xi_\l \xi_j-(E^{(k)}+iF^{(k)})^2|\\
&\le|\sum_{1\le\l, j\le n-1}\frac{a_{\l j}^{(k)}}{a_{nn}^{(k)}}\xi_\l \xi_j|+|(E^{(k)}+iF^{(k)})^2|.
\end{aligned}
\end{equation}
By \eqref{elliptic1}, \eqref{elliptic2}, and \eqref{04-8}, we can obtain
\[
\begin{aligned}
&|\sum_{1\le\l, j\le n-1}\frac{a_{\l j}^{(k)}}{a_{nn}^{(k)}}\xi_\l \xi_j|^2=|\frac{1}{a_{nn}^{(k)}}{\cal A}^{(k)}\xi\cdot\xi|^2\;\;(\mbox{with}\;\xi=(\xi',0))\\
&=|a_{nn}^{(k)}|^{-2}|\cos\theta M^{(k)}\xi\cdot\xi+\gamma \sin\theta N^{(k)}\xi\cdot\xi+i(-\sin\theta M^{(k)}\xi\cdot\xi+\gamma\cos\theta N^{(k)}\xi\cdot\xi)|^2\\
&=|a_{nn}^{(k)}|^{-2}[(M^{(k)}\xi\cdot\xi)^2+\gamma^2(N^{(k)}\xi\cdot\xi)^2]\\
&\le\frac{\Lambda^2(1+\gamma^2)|\xi|^4}{\lambda_0^2(1+\gamma^2)}=\tilde\lambda_1^{-1}|\xi|^4,
\end{aligned}
\]
where we have used the estimate
\begin{equation}\label{ann}
\lambda_0(1+\gamma^2)^{1/2}\le|a_{nn}^{(k)}|\le\Lambda_0(1+\gamma^2)^{1/2}
\end{equation}
in deriving the inequality above. We thus obtain
\begin{equation}\label{annk}
|\sum_{1\le\l, j\le n-1}\frac{a_{\l j}^{(k)}}{a_{nn}^{(k)}}\xi_\l \xi_j|\le\tilde\lambda_1^{-1/2}|\xi'|^2.
\end{equation}
Furthermore, we can estimate
\begin{equation}\label{055}
\begin{aligned}
|(E^{(k)}+iF^{(k)})^2|&=|(\sum_{1\le j\le n-1}\frac{a_{nj}^{(k)}}{a_{nn}^{(k)}}\xi_j)^2|=|\sum_{1\le j\le n-1}\frac{a_{nj}^{(k)}}{a_{nn}^{(k)}}\xi_j|^2\\
&\le(\sum_{1\le j\le n-1}|\frac{a_{nj}^{(k)}}{a_{nn}^{(k)}}|^2)|\xi'|^2\le\frac{(n-1)\Lambda_0^2(1+\gamma^2)}{\lambda_0^2(1+\gamma^2)}|\xi'|^2\\
&=(n-1)\tilde\lambda_1^{-1}|\xi'|^2.
\end{aligned}
\end{equation}
Substituting \eqref{annk}, \eqref{055} into \eqref{a2} gives
\begin{equation}\label{0551}
(A^{(k)})^2+(B^{(k)})^2\le(\tilde\lambda_1^{-1/2}+(n-1)\tilde\lambda_1^{-1})|\xi'|^2\le n\frac{\Lambda_0^2}{\lambda_0^2}|\xi'|^2.
\end{equation}
It follows from \eqref{a1} and \eqref{0551} that
\begin{equation}\label{0552}
\|\delta^{(2)}v_2+\delta^{(1)}v_1\|^2\le 2\Lambda_0^2n\frac{\Lambda_0^2}{\lambda_0^2}|\xi'|^2.
\end{equation}

Next, we want to estimate $\|\eta^{(2)}u_2+\eta^{(1)}u_1\|$ from below. As above, we have
\begin{equation}\label{eta}
\begin{aligned}
\|\eta^{(2)}u_2+\eta^{(1)}u_1\|^2&=(\eta^{(2)})^2[(A^{(2)})^2+(B^{(2)})^2]+(\eta^{(1)})^2[(A^{(1)})^2+(B^{(1)})^2]\\
&+2\eta^{(1)}\eta^{(2)}[(A^{(2)})^2+(B^{(2)})^2]^{1/2}[(A^{(1)})^2+(B^{(1)})^2]^{1/2}\cos\phi.
\end{aligned}
\end{equation}
Recall that $B_1>0$, $B_2<0$. Thus,
\[
\begin{aligned}
\cos\phi&=\frac{A^{(1)}A^{(2)}+B^{(1)}B^{(2)}}{[(A^{(2)})^2+(B^{(2)})^2]^{1/2}[(A^{(1)})^2+(B^{(1)})^2]^{1/2}}\\
&=\frac{A^{(1)}A^{(2)}-|B^{(1)}||B^{(2)}|}{[(A^{(2)})^2+(B^{(2)})^2]^{1/2}[(A^{(1)})^2+(B^{(1)})^2]^{1/2}}\\
&=\frac{1-\frac{|B^{(1)}|}{A^{(1)}}\frac{|B^{(2)}|}{A^{(2)}}}{(1+(\frac{B^{(2)}}{A^{(2)}})^2)^{1/2}(1+(\frac{B^{(1)}}{A^{(1)}})^2)^{1/2}}.
\end{aligned}
\]
Notice that by \eqref{akk} and \eqref{0552}
\[
0\le\frac{|B^{(k)}|}{A^{(k)}}\le \frac{\sqrt{(A^{(k)})^2+(B^{(k)})^2}}{A^{(k)}}\le\frac{\sqrt{n}\frac{\Lambda_0}{\lambda_0}|\xi'|}{\sqrt{\tilde\lambda_1}|\xi'|}=\frac{\sqrt{n}\Lambda_0^2}{\lambda_0^2}:=\tilde\lambda_2\ge 1.
\]
It is readily seen that the function 
\[
f(x,y)=\frac{1-xy}{\sqrt{1+x^2}\sqrt{1+y^2}}
\]
defined on $(x,y)\in[0,\tilde\lambda_2]\times[0,\tilde\lambda_2]$ attains its minimum at $x=y=\tilde\lambda_2$. Hence, we have
\[
\cos\phi\ge\frac{1-\tilde\lambda_2^2}{1+\tilde\lambda_2^2}=-1+\frac{2}{1+\tilde\lambda_2^2}.
\]

Now \eqref{eta} gives
\begin{equation}\label{eta2}
\begin{aligned}
&\|\eta^{(2)}u_2+\eta^{(1)}u_1\|^2\\
\ge&(\eta^{(2)})^2[(A^{(2)})^2+(B^{(2)})^2]+(\eta^{(1)})^2[(A^{(1)})^2+(B^{(1)})^2]\\
&+2\eta^{(1)}\eta^{(2)}[(A^{(2)})^2+(B^{(2)})^2]^{1/2}[(A^{(1)})^2+(B^{(1)})^2]^{1/2}(-1+\frac{2}{1+\tilde\lambda_2^2})\\
=&\left((\eta^{(2)})^2[(A^{(2)})^2+(B^{(2)})^2]^{1/2}-(\eta^{(1)})^2[(A^{(1)})^2+(B^{(1)})^2]^{1/2}\right)^2\\
&+\frac{4}{1+\tilde\lambda_2^2}\eta^{(1)}\eta^{(2)}[(A^{(2)})^2+(B^{(2)})^2]^{1/2}[(A^{(1)})^2+(B^{(1)})^2]^{1/2}\\
&\ge\frac{4}{1+\tilde\lambda_2^2}A^{(1)}A^{(2)}\lambda_0^2\ge \frac{4}{1+\tilde\lambda_2^2}\tilde\lambda_1|\xi'|^2\lambda_0^2=\frac{4}{1+\tilde\lambda_2^2}\frac{\lambda_0^4}{\Lambda_0^2}|\xi'|^2.
\end{aligned}
\end{equation}
Hence, in view of \eqref{0552}, \eqref{eta2}, if we choose
\begin{equation}\label{gamma0}
\gamma_0=\frac{\sqrt{2}\lambda_0^5}{\Lambda_0^3\sqrt{n\lambda_0^4+n^2\Lambda_0^4}},
\end{equation}
then for $\gamma<\gamma_0$ we have
\[
\|\eta^{(2)}u_2+\eta^{(1)}u_1\|^2>\gamma^2\|\delta^{(2)}v_2+\delta^{(1)}v_1\|^2.
\]
In other words, \eqref{tc6} cannot hold (i.e., $\mbox{\rm det}T\ne 0$), and equivalently, \eqref{ne0} is satisfied. In conclusion, we have shown that
\begin{theorem}\label{transthm}
Assume that $a_{\l j}^{(k)}$ have properties \eqref{symm}-\eqref{elliptic2}. Moreover, the number $\gamma$ in \eqref{complex} satisfies $\gamma<\gamma_0$, where $\gamma_0$ is defined in \eqref{gamma0}. Let $\psi_\eps$ be given by \eqref{psieps} with $a_1, a_2$ satisfying \eqref{tc2}. The $\{P_k,\psi_\eps, T_k^j,\; k=1,2,\; j=1,2\}$ satisfies the transmission condition at $0$. 
\end{theorem}

\subsection{Strong pseudoconvexity}\label{stpseudo}

Here we want to check the strong pseudoconvexity condition for the operator ${\mathcal L}_0$ and the weight function $\psi_\eps(x)$ in $B_{\delta'}\cap\Omega_1$ and $B_{\delta'}\cap\Omega_2$ for some small $\delta'>0$. Even though ${\mathcal L}_0$ is represented by $P_k$ in $\Omega_k$, $k=1,2$, it is not necessary to discuss the strong pseudoconvexity condition for $P_1$ and $P_2$ separately.  We suppress the index $k$ in notations and denote the symbol
\[
p(\xi)=\sum_{1\le j,\l\le n}a_{\l j}\xi_\ell \xi_j
\]
with $a_{\l j}=M_{\l j}+i\gamma N_{\l j}$ and consider the weight function
\[
\psi_\eps(x)=\alpha x_n+\frac{\beta}{2}x_n^2-\frac{\eps}{2}|x'|^2.
\]
In view of the definition of $\psi_\eps(x)$ in \eqref{2.1}, $\alpha$ here represents either $\alpha_2=\alpha_+$ or $\alpha_1=\alpha_-$. Hence, we have that
\[
(\partial_j\psi_\eps(x))_{j=1}^n=\nabla\psi_\eps(x)=(-\eps x', \alpha+\beta x_n)
\]
and
\begin{equation}\label{hessian}
(\partial^2_{\l j}\psi_\eps(x))_{\l, j=1}^n=\nabla^2\psi_\eps(x)=\begin{pmatrix}-\eps I_{n-1}&0\\0&\beta\end{pmatrix}.
\end{equation}
The strong pseudoconvexity condition reads that in $B_{\delta'}$, if 
\[
\left\{
\begin{aligned}
&p(\xi+i\tau\nabla\psi_\eps(x))=0,\\
&(\xi,\tau)\ne 0,\;\;\nabla\psi_\eps(x)\ne 0,\;\; x\in\overline{B_{\delta'}},
\end{aligned}\right.
\]
then
\begin{equation}\label{sp}
\begin{aligned}
Q(x,\xi,\tau):=&\sum_{\l, j=1}^n\partial^2_{\l j}\psi_\eps(x)\partial_{\xi_j} p(\xi+i\tau\nabla\psi_\eps(x))\overline{\partial_{\xi_\l}p(\xi+i\tau\nabla\psi_\eps(x))}\\
&+\frac{1}{\tau}\mbox{\rm Im}\sum_{j=1}^n\partial_jp(\xi+i\tau\nabla\psi_\eps(x))\overline{\partial_{\xi_j}p(\xi+i\tau\nabla\psi_\eps(x))}\\
=&\sum_{\l, j=1}^n\partial^2_{\l j}\psi_\eps(x)\partial_{\xi_j} p(\xi+i\tau\nabla\psi_\eps(x))\overline{\partial_{\xi_\l}p(\xi+i\tau\nabla\psi_\eps(x))}>0
\end{aligned}
\end{equation}
(see \cite{ho0}). 

We now write
\[
\begin{aligned}
&p(\xi+i\tau\nabla\psi_\eps)\\
=&\sum_{1\le\l, j\le n}a_{\l j}(\xi_\l+i\tau\partial_\l\psi_\eps)(\xi_\l+i\tau\partial_\l\psi_\eps)\\
=&\sum_{1\le\l, j\le n}a_{\l j}\xi_\l\xi_j+2i\sum_{1\le\l, j\le n}a_{\l j}\xi_l(\tau\partial_j\psi_\eps)-\sum_{1\le \l j\le n}a_{\l j}(\tau\partial_\l\psi_\eps)(\tau\partial_j\psi_\eps).
\end{aligned}
\]
Hence $p(\xi+i\tau\nabla\psi_\eps)=0$ implies
\begin{equation}\label{11}
\sum_{1\le \l j\le n}a_{\l j}(\tau\partial_\l\psi_\eps)(\tau\partial_j\psi_\eps)=\sum_{1\le\l, j\le n}a_{\l j}\xi_\l\xi_j+2i\sum_{1\le\l, j\le n}a_{\l j}\xi_l(\tau\partial_j\psi_\eps)
\end{equation}
By \eqref{complex}-\eqref{elliptic2}, we have
\begin{equation*}
\Big{|}\sum_{1\le \l j\le n}a_{\l j}(\tau\partial_\l\psi_\eps)(\tau\partial_j\psi_\eps)\Big{|}\le\sqrt{1+\gamma^2}\Lambda_0|\tau\nabla\psi_\eps|^2.
\end{equation*}
From this estimate, we obtain from \eqref{11} that
\begin{equation}\label{112}
\begin{aligned}
\sqrt{1+\gamma^2}\Lambda_0|\tau\nabla\psi_\eps|^2\ge &\Big{|}\sum_{1\le\l, j\le n}a_{\l j}\xi_\l\xi_j+2i\sum_{1\le\l, j\le n}a_{\l j}\xi_l(\tau\partial_j\psi_\eps)\Big{|}\\
\ge &\sqrt{1+\gamma^2}\lambda_0|\xi|^2-2\sqrt{1+\gamma^2}\Lambda_0|\xi||\tau\nabla\psi_\eps|\\
\ge&\sqrt{1+\gamma^2}\lambda_0|\xi|^2-\frac{\sqrt{1+\gamma^2}\lambda_0}{2}|\xi|^2-2\frac{\sqrt{1+\gamma^2}\Lambda^2_0}{\lambda_0}|\tau\nabla\psi_\eps|^2,
\end{aligned}
\end{equation}
which leads to
\begin{equation}\label{22}
\frac{\lambda_0}{2}|\xi|^2\le (\Lambda_0+\frac{2\Lambda_0^2}{\lambda_0})|\tau\nabla\psi_\eps|^2.
\end{equation}
By \eqref{22} and exchanging the roles of $\xi$ and $\tau\nabla\psi_\eps$ in \eqref{112}, we thus conclude that there exist positive constants $C_1, C_2$, depending on $\lambda_0, \Lambda_0$ such that
\begin{equation}\label{222}
C_1|\xi|\le|\tau\nabla\psi_\eps|\le C_2|\xi|
\end{equation}
whenever $p(\xi+i\tau\nabla\psi_\eps)=0$. 

As in \eqref{pk} and \eqref{bk}, we can write
\[
p(\xi)=a_{nn}[(\xi_n+\sum_{1\le j\le n-1}\frac{a_{nj}}{a_{nn}}\xi_j)^2+b(\xi')],
\]
where
\[
b(\xi')=\frac{1}{a_{nn}^{2}}\sum_{1\le\l, j\le n-1}(a_{\l j}a_{nn}-a_{n\l}a_{nj})\xi_\l\xi_j.
\]
Similar to \eqref{bk0} and \eqref{efk}, we further express
\begin{equation}\label{4-1}
b(\xi')=(A(\xi')-iB(\xi'))^2,
\end{equation}
with $A(\xi')\ge 0$ and 
\begin{equation}\label{4-2}
\sum_{1\le j\le n-1}\frac{a_{nj}}{a_{nn}}\xi_j=E(\xi')+iF(\xi'),
\end{equation}
where $E(\xi'), F(\xi')\in\R$. 

To verify that \eqref{sp} for $x$ near $0$, we first derive an estimate of $Q(0,\xi,\tau)$. At $x=0$, we have $\partial_j\psi_\eps(0)=0$, $1\le j\le n-1$ and $\partial_n\psi_\eps(0)=\alpha$, i.e.,
\[
\xi+i\tau\nabla\psi_\eps(0)=(\xi',\xi_n+i\tau\alpha).
\]  
Thus, we can rewrite
\begin{equation}\label{p0k}
p(\xi+i\tau\nabla\psi_\eps(0))=p(\xi+i\tau\alpha e_n)=a_{nn}(\xi_n-\sigma_1)(\xi_n-\sigma_2),
\end{equation}
where
\begin{equation}\label{s12}
\left\{
\begin{aligned}
&\sigma_1=-E-B-i(\tau\alpha+F+A),\\
&\sigma_2=-E+B-i(\tau\alpha+F-A).
\end{aligned}\right.
\end{equation}
From now on we suppress the dependence of coefficients at $0$ if there is no danger of causing confusion. 

By \eqref{hessian}, we have that
\begin{equation*}
Q(0,\xi,\tau)=-\eps\sum_{1\le j\le n-1}|\partial_{\xi_j}p(\xi+i\tau\alpha e_n)|^2+\beta|\partial_{\xi_n}p(\xi+i\tau\alpha e_n)|^2,
\end{equation*}
where for $1\le j\le n-1$
\[
\partial_{\xi_j}p(\xi+i\tau\alpha e_n)=2\sum_{1\le\l\le n-1}a_{\l j}\xi_\l+a_{nj}(\xi_n+i\tau\alpha)
\]
and
\[
\partial_{\xi_n}p(\xi+i\tau\alpha e_n)=2\sum_{1\le\l\le n-1}a_{\l n}\xi_\l+a_{nn}(\xi_n+i\tau\alpha).
\]
Therefore, we can write
\begin{equation}\label{q1}
\begin{aligned}
Q(0,\xi,\tau)=&-4\eps\sum_{1\le j\le n-1}\Big{|}\sum_{1\le\l\le n-1}a_{\l j}\xi_\l+a_{nj}(\xi_n+i\tau\alpha)\Big{|}^2\\
&+4\beta\Big{|}\sum_{1\le\l\le n-1}a_{\l n}\xi_\l+a_{nn}(\xi_n+i\tau\alpha)\Big{|}^2.
\end{aligned}
\end{equation}
It follows from \eqref{p0k} and \eqref{s12} that $p(\xi+i\tau\alpha e_n)=0$ if and only if 
\begin{equation}\label{eqi}
\xi_n+i\tau\alpha=-E-B-i(F+A)
\end{equation}
or
\begin{equation}\label{eqii}
\xi_n+i\tau\alpha=-E+B-i(F-A).
\end{equation}
Therefore, if $p(\xi+i\tau\alpha e_n)=0$, then the second term in \eqref{q1} can be further simplified as 
\begin{equation}\label{4beta}
\begin{aligned}
&\Big{|}\sum_{1\le\l\le n-1}a_{\l n}\xi_\l+a_{nn}(\xi_n+i\tau\alpha)\Big{|}^2
=|a_{nn}|^2\Big{|}\sum_{1\le\l\le n-1}\frac{a_{\l n}}{a_{nn}}\xi_\l+(\xi_n+i\tau\alpha)\Big{|}^2\\
&=|a_{nn}|^2|E+iF+(\xi_n+i\tau\alpha)|^2=|a_{nn}|^2(A^2+B^2),
\end{aligned}
\end{equation}
where we have used \eqref{4-2}, \eqref{eqi} or \eqref{eqii}. Combining \eqref{4-1}, \eqref{4-2}, \eqref{eqi} or \eqref{eqii}, we have that
\begin{equation}\label{ebfa}
|\xi_n+i\tau\alpha|\le|E|+|B|+|F|+|A|\le C\Lambda_0|\xi'|,
\end{equation}
which implies
\begin{equation}\label{2-20}
\sum_{1\le j\le n-1}\Big{|}\sum_{1\le\l\le n-1}a_{\l j}\xi_\l+a_{nj}(\xi_n+i\tau\alpha)\Big{|}^2\le C\Lambda_0^2|\xi'|^2.
\end{equation}

Putting \eqref{q1}, \eqref{4beta}, and \eqref{2-20} together gives
\begin{equation}\label{q11}
Q(0,\xi,\tau)\ge 4\beta|a_{nn}|^2(A^2+B^2)-4\eps C\Lambda_0|\xi'|^2.
\end{equation}
Recall the estimate \eqref{akk} in Lemma~\ref{ak}
\[
A^2\ge\tilde\lambda_1|\xi'|^2+|F|^2\ge\tilde\lambda_1|\xi'|^2.
\]
Using this estimate in \eqref{q11} and choosing $\eps$ sufficiently small leads to 
\begin{equation*}
Q(0,\xi,\tau)\ge 4(\beta\tilde\lambda_1\lambda_0^2-\eps C\Lambda_0)|\xi'|^2\ge 2\beta\tilde\lambda_1\lambda_0^2|\xi'|^2,
\end{equation*}
whenever $p(\xi+i\tau\alpha e_n)=0$. Furthermore, \eqref{ebfa} implies
\[
|\xi+i\tau\alpha e_n|^2\le (1+C^2\Lambda_0^2)|\xi'|^2,
\]
and it follows that if $p(\xi+i\tau\alpha e_n)=0$ then
\begin{equation}\label{XI-2}
Q(0,\xi,\tau)\ge C\beta|\xi+i\tau\alpha e_n|^2.
\end{equation}
In conclusion, we have shown that 
\begin{equation}\label{14-2}
(\xi,\tau)\in\{(\xi,\tau)\in S: p(\xi+i\tau\alpha e_n)=0\}\Rightarrow Q(0,\xi,\tau)>0,
\end{equation}
where $S:=\{(\xi,\tau)\in\R^{n+1}: |\xi|^2+\tau^2=1\}$. 

Now we recall the following elementary theorem. Let $X$ be a compact subset of $\R^N$ and $F, G: X\to\R$ be two continuous functions, then the following two statements are equivalent:
\begin{description}
\item[(i)] $F(x)=0,\;\forall\; x\in X\Rightarrow G(x)>0$.
\item[(ii)] There exist positive constants $C_1, C_2$ such that $C_1G(x)+|F(x)|\ge C_2$, $\forall\; x\in X$.
\end{description}
With the help of this theorem, \eqref{14-2} is equivalent to 
\begin{equation}\label{14-3}
C_1Q(0,\xi,\tau)+|p(\xi+i\tau\alpha e_n)|\ge C_2
\end{equation}
for all $(\xi,\tau)\in S$. Thanks to \eqref{14-3}, we can estimate
\[
\begin{aligned}
C_1Q(x,\xi,\tau)+|p(\xi+i\tau\nabla\psi_\eps(x))|&=C_1Q(0,\xi,\tau)+|p(\xi+i\tau\alpha e_n)|+R(x,\xi,\tau)\\
&\ge C_2+R(x,\xi,\tau),
\end{aligned}
\]
where
\[
R(x,\xi,\tau)=C_1[Q(x,\xi,\tau)-Q(0,\xi,\tau)]+|p(\xi+i\tau\nabla\psi_\eps(x))|-|p(\xi+i\tau\alpha e_n)|.
\]
Observe that $R(0,\xi,\tau)=0$ for $(\xi,\tau)\in S$. Since $R$ is continuous, there exists a small number $\delta'>0$ such that
\[
|R(x,\xi,\tau)|\le \frac{C_2}{2}
\] 
for all $x$ with $|x|\le\delta'<\frac{\alpha}{2\beta}$ and $(\xi,\tau)\in S$. In other words, we have that
\begin{equation}\label{14-5}
C_1Q(x,\xi,\tau)+|p(\xi+i\tau\nabla\psi_\eps(x))|\ge\frac{C_2}{2}
\end{equation}
in $\{|x|\le\delta'\}\times S$. By the elementary theorem stated above, \eqref{14-5} is equivalent to 
\[
\begin{aligned}
&p(\xi+i\tau\nabla\psi_\eps(x))=0,\;\forall\; x\in\overline{B_{\delta'}},\; (\xi,\tau)\in S\\
\Rightarrow&\; Q(x,\xi,\tau)>0,\;\forall\; x\in\overline{B_{\delta'}},\; (\xi,\tau)\in S,
\end{aligned}
\]
which immediately implies the strong pseudoconvexity condition near $0$ in view of the homogeneity of $p$ and $Q$ in $(\xi,\tau)$. 

Having verified the strong psudoconvexity in a neighborhood of $0$ and the transmission conditions at $0$, we can derive a Carleman estimate with weight $\psi_{\varepsilon}(x)$ for the operator ${\mathcal{L}_0}$.
\begin{theorem}{\rm\cite[Theorem 1.6]{BL}} \label{coeffcost}
Assume that coefficients $A_{\pm}(0)$ satisfy conditions \eqref{symm}-\eqref{elliptic2}. There exist $\alpha_+,\alpha_-,\beta,\eps_0, \gamma_0, r_0$ and $C$, depending on $\lambda_0, \Lambda_0$, such that if  $\eps\le\eps_0$, $\gamma\le\gamma_0$, $\tau\geq C$, then
\begin{eqnarray}\label{8.240}
&&\sum_{\pm}\sum_{k=0}^2\tau^{3-2k}\int_{\mathbb{R}^n_{\pm}}|D^k{u}_{\pm}|^2e^{2\tau\psi_{\varepsilon,\pm}(x)}dx+\sum_{\pm}\sum_{k=0}^1\tau^{3-2k}\int_{\mathbb{R}^{n-1}}|D^k{u}_{\pm}(x',0)|^2e^{2\psi_{\varepsilon}(x',0)}dx'\nonumber\\
&&+\sum_{\pm}\tau^2[e^{\tau\psi_{\varepsilon}(\cdot,0)}u_{\pm}(\cdot,0)]^2_{1/2,\mathbb{R}^{n-1}}+\sum_{\pm}[D(e^{\tau\psi_{\varepsilon,\pm}}u_{\pm})(\cdot,0)]^2_{1/2,\mathbb{R}^{n-1}}\nonumber\\
&&\leq C\left(\int_{\mathbb{R}_{\pm}^n}|\mathcal{L}_0(D)(u_{\pm})|^2\,e^{2\tau\psi_{\varepsilon,\pm}(x)}dx+[e^{\tau\psi_{\varepsilon}(\cdot,0)}h^{(0)}_1]^2_{1/2,\mathbb{R}^{n-1}}\right.\\
&&\left.+[D_{x'}(e^{\tau\psi_{\varepsilon}}h^{(0)}_0)(\cdot,0)]^2_{1/2,\mathbb{R}^{n-1}}+\tau^{3}\int_{\mathbb{R}^{n-1}}|h^{(0)}_0|^2e^{2\tau\psi_{\varepsilon}(x,0)}dx+\tau\int_{\mathbb{R}^{n-1}}|h^{(0)}_1|^2e^{2\tau\psi_{\varepsilon}(x,0)}dx\right).\nonumber
\end{eqnarray}
for $u=H_+u_++H_-u_-$,  $u_{\pm}\in C^\infty(\mathbb{R}^{n})$ and ${\rm supp}\, u\subset B'_{r_0}\times[-r_0, r_0]$, and
\begin{equation*}
h^{(0)}_0(x'):=u_+(x',0)-u_-(x',0),\ \forall\, x'\in \mathbb{R}^{n-1},
\end{equation*}
\begin{equation*}
h^{(0)}_1(x'):=A_+(0)\nabla u_+(x',0)\cdot e_n-A_-(0)\nabla u_-(x',0)\cdot e_n,\ \forall\, x\in \mathbb{R}^{n-1}.
\end{equation*}
\end{theorem}

\section{Derivation of the Carleman estimate}\label{sec4}

This section is devoted to the derivation of the Carleman estimate \eqref{8.24} following the ideas used in \cite{dflvw}. We first introduce the partition of unity given in \cite{dflvw}.  For any $r>0$ and $x'\in \R^{n-1}$, denote the $(n-1)$-cube
$Q_{r}(x')=\{y'\in\mathbb{R}^{n-1}:|y'_j-x'_j|\leq r,\,j=1,2,\cdots ,n-1\}$. Let $\vartheta_0\in C^\infty_0(\mathbb{R})$ such that
\begin{equation}\label{tetazero}
0\leq\vartheta_0\leq 1,\quad  {\rm supp}\,\vartheta_0\subset (-3/2,3/2)\,\,\mbox{ and }\,\,\vartheta_0(t)=1\mbox{ for }t\in[-1,1].
\end{equation}
Let  $\vartheta(x')=\vartheta_0(x_1)\cdots\vartheta_0(x_{n-1})$, so that
\begin{equation*}
{\rm supp}\,\vartheta\subset \stackrel{\circ}{Q}_{3/2}(0)\,\,\mbox{ and }\,\,\vartheta(x')=1\mbox{ for }x'\in Q_1(0),
\end{equation*}
where $\stackrel{\circ}{Q}$ denotes the interior of the set $Q$. Given $\mu\geq 1$ and $g\in \mathbb{Z}^{n-1}$, we define
$$x'_g=\frac{g}{\mu}$$
and
$$\vartheta_{g,\mu}(x')=\vartheta(\mu(x'-x'_g)).$$
Thus, we can see that
\begin{equation*}
{\rm supp}\,\vartheta_{g,\mu}\subset \stackrel{\circ}{Q}_{3/2\mu}(x'_g)\subset Q_{2/\mu}(x'_g)
\end{equation*}
and
\begin{equation}\label{6.4}
|D^k\vartheta_{g,\mu}|\leq C_1\mu^k(\chi_{ Q_{3/2\mu}(x'_g)}-\chi_{Q_{1/\mu}(x'_g)}),\quad k=0,1,2,
\end{equation}
where $C_1\geq 1$ depends only on $n$.

Notice that, for any $g\in \mathbb{Z}^{n-1}$,
\begin{equation}\label{cardAg}
card\left(\{g'\in \mathbb{Z}^{n-1}\,:\,{\rm supp}\,\vartheta_{g',\mu}\cap{\rm supp}\,\vartheta_{g,\mu}\neq\emptyset\}\right)=5^{n-1}.\end{equation}
 Thus, we can define
\begin{equation}\label{6.5}
\bar{\vartheta}_{\mu}(x'):=\sum_{g\in \mathbb{Z}^{n-1}}\vartheta_{g,\mu}\geq 1,\quad x'\in \mathbb{R}^{n-1}.
\end{equation}
By \eqref{6.4}, we get that
\begin{equation}\label{6.6}
|D^k\bar{\vartheta}_{\mu}|\leq C_2\mu^k,
\end{equation}
where $C_2\geq 1$ depends on $n$. Define
\begin{equation}\label{etagmu}
\eta_{g,\mu}(x')=\vartheta_{g,\mu}(x')/\bar{\vartheta}_{\mu}(x'),\quad x'\in \mathbb{R}^{n-1},\end{equation}
then we have that
\begin{equation}\label{6.7}
\begin{cases}
\sum_{g\in \mathbb{Z}^{n-1}}\eta_{g,\mu}= 1,\quad x'\in \mathbb{R}^{n-1},\\
{\rm supp}\,\eta_{g,\mu}\subset Q_{3/2\mu}(x'_g)\subset Q_{2/\mu}(x'_g),\\
|D^k\eta_{g,\mu}|\leq C_3\mu^k\chi_{Q_{3/2\mu}(x'_g)},\quad k=0,1,2,
\end{cases}
\end{equation}
where $C_3\geq 1$ depends on $n$.

We will first extend \eqref{8.240} to operators with leading coefficients depending on the vertical variable $x_n$. To do so, we need to derive an interior Carleman estimate for second order elliptic operators having Lipschitz leading coefficients and with the weight function $\psi_\eps$. 
To derive such Carleman estimate, we define the $n$-cube $K_R=\{x=(x_1,\cdots,x_n): |x_j|\le R, 1\le j\le n\}$ for $R>0$. Let us denote 
\[
P(x,D)=\sum_{1\le j, \l\le n}a_{j\l}(x)D_{j\l}^2
\]
and its symbol $p(x,\xi)=\sum_{1\le j, \l\le n}a_{j\l}(x)\xi_j\xi_\l$. Assume that for all $1\le j, \l\le n$ and $x, y\in K_1$,
\begin{equation}\label{cond4}
\left\{\begin{aligned}
&a_{j\l}(x)=a_{\l j}(x),\\
&|a_{j\l}(x)|\le\Lambda,\\
&|a_{j\l}(x)-a_{j\l}(y)|\le M_0|x-y|,\\
&|p(x,\xi)|\ge\lambda|\xi|^2,\;\;\forall\;\xi\in\R^n,
\end{aligned}\right.
\end{equation}
where $\Lambda, \lambda>0$. Let $\varphi(x)\in C^2(K_1)$ be real-valued and satisfy $|\nabla\varphi(x)|\ne 0$ for all $x\in K_1$. We denote
\[
S(x,y;\xi,\tau)=\sum_{\l, j=1}^n\partial^2_{\l j}\varphi(x)\partial_{\xi_j} p(y,\xi+i\tau\nabla\varphi(x))\overline{\partial_{\xi_\l}p(y,\xi+i\tau\nabla\varphi(x))}
\] 
for $x, y\in K_1$, $\xi\in\R^n$, $\tau>0$. 
\begin{pr}\label{pr1}
Assume that the following condition holds:
\begin{equation}\label{stcon}
\left.\begin{aligned}
&p(0,\xi+i\tau\varphi(0))=0\\
&(\xi,\tau)\ne(0,0)\end{aligned}\right\}\Rightarrow S(0,0;\xi,\tau)>0.
\end{equation}
Then there exist $\bar{R}\in(0,1]$, $\delta_0\in(0,1]$, $C_0\ge 1$, $\tau_0\ge 1$, depending on $\lambda, \Lambda, M_0, \|\varphi\|_{C^2(Q_1)}$,  such that
\begin{equation}\label{car00}
\sum_{|\alpha|\le 2}\tau^{3-2|\alpha|}\int|D^\alpha u|^2e^{2\tau\varphi(x)}dx\le C_0\int|P(\delta x,D)u|^2e^{2\tau\varphi(x)}dx,
\end{equation}
$\forall\, u\in C_0^\infty(\stackrel{\circ}{K}_{\bar R})$, $\tau\ge\tau_0$, $0<\delta\le\delta_0$.  
\end{pr}
\pf. In view of the homogeneity in $(\xi,\tau)$, \eqref{stcon} is equivalent to that there exist $C_1>0, C_2>0$ such that
\begin{equation*}
C_2|p(0,\xi+i\tau\nabla\varphi(0)|^2+(|\xi|^2+\tau^2)S(0,0;\xi,\tau)\ge C_1(|\xi|^2+\tau^2)^2,\;\forall\;(\xi,\tau)\in \R^{n+1}.
\end{equation*}
From \eqref{cond4}, we can see that there exists $\bar R\in(0,1]$ such that  
\begin{equation}\label{stcon2}
\tilde C_2|p(y,\xi+i\tau\nabla\varphi(x)|^2+(|\xi|^2+\tau^2)S(x,y;\xi,\tau)\ge \tilde C_1(|\xi|^2+\tau^2)^2,\;\forall\;x, y\in K_{\bar R},\,\forall\;(\xi,\tau)\in \R^{n+1},
\end{equation}
where $\tilde C_1>0$, $\tilde C_2>0$ are independent of $x$, $y$. Thanks to \eqref{stcon2}, the Carleman derived in \cite[Theorem~8.3.1]{ho0} holds for 
\[
P(\delta y,D_x)u=\sum_{1\le\l, j\le n}a_{j\l}(\delta y)D^2_{x_jx_\l}u(x),
\]
that is,
\begin{equation}\label{car22}
\sum_{|\alpha|\le 2}\tau^{3-2|\alpha|}\int|D_x^\alpha u|^2e^{2\tau\varphi(x)}dx\le C_3\int|P(\delta y,D_x)u|^2e^{2\tau\varphi(x)}dx
\end{equation}
for all $u\in C_0^\infty(\overset{\circ}{K}_{\bar R})$, $0\le\delta\le 1$, and $\tau\ge\tau_1$, where $C_3$ and $\tau_1$ do not depend on $\delta$ and $y$. Note that for fixed $\delta$, $y$, $P(\delta y,D_x)$ is an operator having constant coefficients. 

Now we use the partition of unity introduced above, but with $n-1$ being replaced by $n$. In particular, for $h\in{\mathbb Z}^n$, we define
\[
x_h=\frac{h}{\mu},\quad\mu=\sqrt{\eps\tau}\;\;\mbox{with}\;\;\tau\ge\frac{1}{\eps},
\]
where $\eps\in(0,1]$ will be chosen later. Let $u\in C_0^\infty(\stackrel{\circ}{K}_{\bar R})$, in view of the first relation in \eqref{6.7}, we have
\[
u(x)=\sum_{h\in{\mathbb Z}^n}u(x)\eta_{h,\mu}(x),
\]
where $\eta_{h,\mu}(x)$ is defined similarly as in \eqref{etagmu} with $n-1, g$ being replaced by $n, h$, respectively. Applying \eqref{car22} with $y=x_h$ implies
\begin{equation}\label{car33}
\begin{aligned}
&\sum_{|\alpha|\le 2}\tau^{3-2|\alpha|}\int|D^\alpha u|^2e^{2\tau\varphi(x)}dx\\
\le &c\sum_{h\in{\mathbb Z}^n}\sum_{|\alpha|\le 2}\tau^{3-2|\alpha|}\int|D^\alpha (u\eta_{h,\mu})|^2e^{2\tau\varphi(x)}dx\\
\le&cC_3\int|P(\delta x_h,D)(u\eta_{h,\mu})|^2e^{2\tau\varphi(x)}dx,\;\;\forall\;\tau\ge\tau_2=\min\{\tau_1,\frac{1}{\eps}\},
\end{aligned}
\end{equation}
where $c=c(n)$. 

Now we write
\begin{equation}\label{est1}
|P(\delta x_h,D)(u\eta_{h,\mu})|\le|P(\delta x,D)(u\eta_{h,\mu})|+|(P(\delta x_h,D)-P(\delta x,D))(u\eta_{h,\mu})|
\end{equation}
and use \eqref{6.6}, the second inequality of \eqref{cond4}, to estimate
\begin{equation}\label{est2}
|P(\delta x,D)(u\eta_{h,\mu})|\le|P(\delta x,D)u|\eta_{h,\mu}+C_4\Lambda(\sqrt{\eps\tau}|Du|+\eps\tau|u|)\chi_{K_{2/\mu}(x_h)}
\end{equation}
and
\begin{equation}\label{est3}
\begin{aligned}
&|(P(\delta x_h,D)-P(\delta x,D))(u\eta_{h,\mu})|=|\sum_{1\le j\l\le n}(a_{j\l}(\delta x_h)-a_{j\l}(\delta x))D^2_{j\l}(u\eta_{h,\mu})|\\
\le&\,\eta_{h,\mu}\sum_{1\le j\l\le n}|a_{j\l}(\delta x_h)-a_{j\l}(\delta x)||D^2_{j\l} u|+2C_4\Lambda(\sqrt{\eps\tau}|Du|+\eps\tau|u|)\chi_{K_{2/\mu}(x_h)}\\
\le &c\,\eta_{h,\mu}\frac{\delta M_0}{\mu}|D^2u|+2C_4\Lambda(\sqrt{\eps\tau}|Du|+\eps\tau|u|)\chi_{K_{2/\mu}(x_h)}
\end{aligned}
\end{equation}
with $c=c(n)$. Here $K_{2/\mu}(x_h)$  denotes the $n$-cube centered at $x_h$ with length $4/\mu$ and $\chi_{K_{2/\mu}(x_h)}$ is the characteristic function of $K_{2/\mu}(x_h)$. Substituting \eqref{est1}-\eqref{est3} into \eqref{car33} gives
\begin{equation}\label{car55}
\begin{aligned}
&\sum_{|\alpha|\le 2}\tau^{3-2|\alpha|}\int|D^\alpha u|^2e^{2\tau\varphi(x)}dx\\
\le&C_5\int|P(\delta x,D)u|^2e^{2\tau\varphi(x)}dx\\
&+C_5\left\{\frac{\delta^2M_0^2}{\eps\tau}\int|D^2u|^2e^{2\tau\varphi(x)}dx+\eps\tau\int|Du|^2e^{2\tau\varphi(x)}dx+(\eps\tau)^2\int|u|^2e^{2\tau\varphi(x)}dx\right\}
\end{aligned}
\end{equation}
for all $\tau\ge\tau_2$, where $C_5\ge 1$. Finally, by choosing $\eps=1/(2C_5)$ and $\delta_0=\eps$, all terms inside of the curved brace on the right hand side of \eqref{car55} can be absorbed by its left hand side and \eqref{car00} follows immediately. \eproof

\subsection{Carleman estimate for operators depending on the vertical variable}

Here we would like to prove a Carleman estimate for the operator that satisfies conditions \eqref{symm}-\eqref{elliptic2} but depending only on the $x_n$ variable. That is, we consider
\begin{equation*}
\mathcal{L}(x_n,D)u:=\sum_{\pm}H_{\pm}{\rm div}(A_{\pm}(x_n)\nabla u_{\pm}),
\end{equation*}
where $u_{\pm}\in C^\infty(\mathbb{R}^{n})$ and ${\rm supp}\, u\subset B'_{r_0}\times[-r_0, r_0]$, where $r_0$ is the number obtained in Theorem~\ref{coeffcost}. Introduce $\delta\in (0,1)$ that will be chosen later, define
\begin{equation*}
\phi_{\delta}(x):=\psi_{\delta}(\delta^{-1}x)=\psi_{\delta}(\delta^{-1}x',\delta^{-1}x_n),
\end{equation*}
and consider the scaled operator
$$\mathcal{L}(\delta x_n,D)u:=
\sum_{\pm}H_{\pm}{\rm div}(A_{\pm}(\delta x_n)\nabla u_{\pm}).$$ Notice that $A_\pm(\delta x_n)$ satisfies assumptions \eqref{elliptic1}, \eqref{elliptic2} and also the Lipschitz condition
\begin{equation}\label{lipdelta}
|A_{\pm}(\delta \tilde x_n)-A_{\pm}(\delta x_n)|\leq M_0\delta|\tilde x_n-x_n|.
\end{equation}

Let $\vartheta_0\in C^\infty_0(\mathbb{R})$ be given as in \eqref{tetazero}. For $\mu\geq 1$ satisfying ${2}/{\mu}<r_0$,  we define
\begin{equation}\label{etamu}
\eta_{\mu}(x_n)=\vartheta_0(\mu x_n),
\end{equation}
\begin{equation}\label{v-z}
v_{\mu}(x',x_n)=\eta_{\mu}(x_n)u(x',x_n),\,\, \mbox{ and }\, z_{\mu}(x',x_n)=(1-\eta_{\mu}(x_n))u(x',x_n).
\end{equation}
Since $v_{\mu,\pm}(x',0)=u_{\pm}(x',0)$ and $\nabla v_{\mu,\pm}(x',0)=\nabla u_{\pm}(x',0)$, we have trivially 
\begin{equation}\label{7.5-v}
v_{\mu,+}(x',0)-v_{\mu,-}(x',0)=u_+(x',0)-u_-(x',0)=h^{(0)}_0(x'),\ \forall\, x'\in \mathbb{R}^{n-1}
\end{equation}
and
\begin{equation}\label{7.6-v}
\begin{aligned}
&A_+(0)\nabla v_{\mu,+}(x',0)\cdot e_n-A_-(0)\nabla v_{\mu,-}(x',0)\cdot e_n\\
=& A_+(0)\nabla u_{+}(x',0)\cdot e_n-A_-(0)\nabla u_{-}(x',0)\cdot e_n=h^{(0)}_1(x'),\; \forall\, x'\in \mathbb{R}^{n-1}.
\end{aligned}
\end{equation}
The aim of this section is to prove a simple version of \eqref{8.24}:
\begin{eqnarray}\label{8.241}
&&\sum_{\pm}\sum_{k=0}^2\tau^{3-2k}\int_{\mathbb{R}^n_{\pm}}|D^k{u}_{\pm}|^2e^{2\tau\psi_{\varepsilon,\pm}(x)}dx+\sum_{\pm}\sum_{k=0}^1\tau^{3-2k}\int_{\mathbb{R}^{n-1}}|D^k{u}_{\pm}(x',0)|^2e^{2\psi_{\varepsilon}(x',0)}dx'\nonumber\\
&&+\sum_{\pm}\tau^2[e^{\tau\psi_{\varepsilon}(\cdot,0)}u_{\pm}(\cdot,0)]^2_{1/2,\mathbb{R}^{n-1}}+\sum_{\pm}[D(e^{\tau\psi_{\varepsilon,\pm}}u_{\pm})(\cdot,0)]^2_{1/2,\mathbb{R}^{n-1}}\nonumber\\
&&\leq C\left(\int_{\mathbb{R}_{\pm}^n}|\mathcal{L}(\delta x_n,D)(u_{\pm})|^2\,e^{2\tau\psi_{\varepsilon,\pm}(x)}dx+[e^{\tau\psi_{\varepsilon}(\cdot,0)}h^{(0)}_1]^2_{1/2,\mathbb{R}^{n-1}}\right.\\
&&\left.+[D_{x'}(e^{\tau\psi_{\varepsilon}}h^{(0)}_0)(\cdot,0)]^2_{1/2,\mathbb{R}^{n-1}}+\tau^{3}\int_{\mathbb{R}^{n-1}}|h^{(0)}_0|^2e^{2\tau\psi_{\varepsilon}(x,0)}dx+\tau\int_{\mathbb{R}^{n-1}}|h^{(0)}_1|^2e^{2\tau\psi_{\varepsilon}(x,0)}dx\right).\nonumber
\end{eqnarray}

To proceed the proof of \eqref{8.241}, we first note that ${\rm supp}\, z_\mu\subset B'_{r_0}\times[-r_0, r_0]$ and vanishes in the strip $\mathbb{R}^{n-1}\times [-\frac{1}{\mu},\frac{1}{\mu}]$. It is clear that $A_\pm(\delta x_n)$ satisfies \eqref{elliptic1}, \eqref{elliptic2} and \eqref{7.4}. Moreover, estimate \eqref{14-3} implies that the condition \eqref{stcon} holds for $\sum_{1\le j, \l\le n}a^{\pm}_{j\l}(x)D_{j\l}^2$ with $\varphi=\psi_\eps$. Observe that $z_\mu$ is supported away from $x_n=0$. Therefore, it follows from \eqref{car00} in Proposition~\ref{pr1} that there exist $\delta_0\in(0,1]$, $\tau_0>0$, and choose a small $r_0$ if necessary, such that   
\begin{equation}\label{8.24interno}
\sum_{k=0}^2\tau^{3-2k}\int_{\mathbb{R}^n}|D^kz_{\mu}|^2e^{2\tau\psi_{\varepsilon}(x)}dx
\leq 
C\int_{\mathbb{R}^n}|\mathcal{L}(\delta x_n,D)z_{\mu}|^2\,e^{2\tau\psi_{\varepsilon}(x)}dx
\end{equation}
for all $\tau\ge\tau_0$, $0<\delta\le\delta_0$, where $C$ depends on $\Lambda_0$, $\lambda_0$, and $M_0$.

Let us denote by $LHS(u)$ the left hand side of inequality \eqref{8.241}. We have
\begin{equation}\label{6-178c}
\begin{aligned}
LHS(u)&\leq 2\left(LHS(v_\mu)+LHS(z_\mu)\right)\\&=2\left(LHS(v_\mu)+\sum_{k=0}^2\tau^{3-2k}\int_{\mathbb{R}^n}|D^kz_{\mu}|^2e^{2\tau\psi_{\varepsilon}(x)}dx\right).
\end{aligned}
\end{equation}
Then applying \eqref{8.240} to $v_\mu$ and using \eqref{8.24interno} leads to
\begin{equation}\label{6-178cdue}
\begin{aligned}
LHS(u)\leq &C\left(\int_{\mathbb{R}^n}|\mathcal{L}_0(D)v_\mu|^2\,e^{2\tau\psi_{\varepsilon}(x)}dx+[e^{\tau\psi_{\varepsilon}(\cdot,0)}h^{(0)}_1]^2_{1/2,\mathbb{R}^{n-1}}\right.\\
&+[D_{x'}(e^{\tau\psi_{\varepsilon}}h^{(0)}_0)(\cdot,0)]^2_{1/2,\mathbb{R}^{n-1}}+\tau^{3}\int_{\mathbb{R}^{n-1}}|h^{(0)}_0|^2e^{2\tau\psi_{\varepsilon}(x',0)}dx'\\&\left.+\tau\int_{\mathbb{R}^{n-1}}|h^{(0)}_1|^2e^{2\tau\psi_{\varepsilon}(x',0)}dx'+\int_{\mathbb{R}^n}|\mathcal{L}(\delta x_n,D)z_{\mu}|^2\,e^{2\tau\psi_{\varepsilon}(x)}dx\right).
\end{aligned}
\end{equation}
By \eqref{elliptic1}, \eqref{elliptic2}, \eqref{lipdelta} and \eqref{etamu} and since $\mu>1$, we can estimate
\begin{eqnarray}\label{1-179c}
&&\left|\mathcal{L}_0(D)v_\mu\right|\nonumber\\
&\leq& \left|\mathcal{L}(\delta x_n,D)v_\mu\right|+\left|\mathcal{L}(\delta x_n,D)v_\mu-\mathcal{L}_0(D)v_\mu\right|\nonumber\\
 &\leq&\left|\mathcal{L}(\delta x_n,D)u\right|\eta_\mu+\frac{2\delta M_0}{\mu}\sum_{\pm}\left|D^2u_\pm\right|\eta_\mu\nonumber\\&&
+C(\delta M_0+\Lambda_0)\sum_\pm(\mu|D u_\pm|+\mu^2|u_\pm|)\chi_{\phantom{l}_{\mathbb{R}^{n-1}\times\left( \left[-\frac{2}{\mu},\frac{2}{\mu}\right]\setminus\left[-\frac{1}{\mu},\frac{1}{\mu}\right]\right)}}.
\end{eqnarray}
On the other hand, we have
\begin{eqnarray}\label{1-180}
&&\left|\mathcal{L}(\delta x_n,D)z_\mu\right|\nonumber\\
&\leq& \left|\mathcal{L}(\delta x_n,D)u\right|(1-\eta_\mu)\nonumber\\
&&+C(\delta M_0+\Lambda_0)\sum_\pm\left(\mu|D u_\pm|+\mu^2|u_\pm|\right)\chi_{\phantom{l}_{\mathbb{R}^{n-1}\times \left( \left[-\frac{2}{\mu},\frac{2}{\mu}\right]\setminus\left[-\frac{1}{\mu},\frac{1}{\mu}\right]\right)}}.
\end{eqnarray}
Putting \eqref{1-179c}, \eqref{1-180}, and \eqref{6-178cdue} together implies
\begin{equation}\label{sn}
LHS(u)\leq C_1\left(\int_{\mathbb{R}^{n}} \left|\mathcal{L}(\delta x_n,D)u\right|^2e^{2\tau\psi_{\varepsilon}(x)}dx+\mathcal{T}_R\right)+C_2\mathcal{R}
\end{equation}
where
\begin{eqnarray*}
	 \mathcal{T}_R&=&[e^{\tau\psi_{\varepsilon}(\cdot,0)}h^{(0)}_1]^2_{1/2,\mathbb{R}^{n-1}}+[D_{x'}(e^{\tau\psi_{\varepsilon}}h^{(0)}_0)(\cdot,0)]^2_{1/2,\mathbb{R}^{n-1}}\\&&+\tau^{3}\int_{\mathbb{R}^{n-1}}|h^{(0)}_0|^2e^{2\tau\psi_{\varepsilon}(x',0)}dx'+\tau\int_{\mathbb{R}^{n-1}}|h^{(0)}_1|^2e^{2\tau\psi_{\varepsilon}(x',0)}dx',
\end{eqnarray*}
\begin{equation*}
\mathcal{R}=\frac{\delta^2}{\mu^2}\sum_\pm\int_{\mathbb{R}^n_{\pm}}|D^2{u}_{\pm}|^2dx+\mu^2\sum_\pm\int_{\mathbb{R}^n_{\pm}}|D u_\pm|^2dx+\mu^4\int_{\mathbb{R}^{n}}|u|^2dx,
\end{equation*}
$C_1$ depends only on $\Lambda_0$ and $\lambda_0$ and $C_2$ depends only on $\Lambda_0$, $\lambda_0$ and $M_0$.

Now we choose $\mu=\sqrt{ \eps\tau}$ and calculate
\begin{eqnarray}\label{ass}
LHS(u)-C_2\mathcal{R}&=&\frac{1}{\tau}\left(1-\frac{C_2\delta^2}{\eps}\right)\sum_{\pm}\int_{\mathbb{R}^n_{\pm}}|D^2{u}_{\pm}|^2e^{2\tau\psi_{\eps}}dx\nonumber\\
&&+\tau\left(1-C_2\eps\right)\sum_{\pm}\int_{\mathbb{R}^n_{\pm}}|D{u}_{\pm}|^2e^{2\tau\psi_{\varepsilon}}dx\nonumber\\&&+\tau^3\left(1-\frac{C_2\eps}{\tau}\right)\int_{\mathbb{R}^n}|u|^2e^{2\tau\psi_{\varepsilon}}dx+\mathcal{T}_L,
\end{eqnarray}
where
\begin{eqnarray*}
\mathcal{T}_L&=&\sum_{\pm}\sum_{k=0}^1\tau^{3-2k}\int_{\mathbb{R}^{n-1}}|D^k{u}_{\pm}(x',0)|^2e^{2\psi_{\varepsilon}(x',0)}dx'\nonumber\\
&&+\sum_{\pm}\tau^2[e^{\tau\psi_{\varepsilon}(\cdot,0)}u_{\pm}(\cdot,0)]^2_{1/2,\mathbb{R}^{n-1}}+\sum_{\pm}[D(e^{\tau\psi_{\varepsilon,\pm}}u_{\pm})(\cdot,0)]^2_{1/2,\mathbb{R}^{n-1}}.
\end{eqnarray*}
By choosing $\eps$ and $\delta$ satisfying
\begin{equation}\label{deltaep}
\delta^2\leq \frac{\eps}{2C_2}\quad\mbox{ and }\quad\eps\leq  \frac{1}{2C_2},
\end{equation}
estimate \eqref{8.241} follows easily from \eqref{sn} and \eqref{ass}.

%
%
%

\subsection{Carleman estimate for operators depending on all variables}

We now want to extend the estimate \eqref{8.241} to operators with coefficients depending also on the variables $x'$. To treat this case we proceed exactly as in \cite[Section 4.2, pp.198-200]{dflvw}, that is, we approximate with coefficients depending only on $x_n$. We use the partition of unity introduced at the beginning of Section~\ref{sec4} and show that
\begin{equation}\label{4.25}
LHS(u)\leq C \sum_{g\in \mathbb{Z}^{n-1}} LHS\left(u\eta_{g,\mu}\right)+CR_1,
\end{equation}
where we define
\[
\begin{aligned}
LHS(u)=&\sum_{\pm}\sum_{k=0}^2\tau^{3-2k}\int_{\mathbb{R}^n_{\pm}}|D^k{u}_{\pm}|^2e^{2\tau\psi_{\eps,\pm}(x',x_n)}dx'dx_n\\
&+\sum_{\pm}\sum_{k=0}^1\tau^{3-2k}\int_{\mathbb{R}^{n-1}}|D^k{u}_{\pm}(x',0)|^2e^{2\psi_\eps(x',0)}dx'\\
&+\sum_{\pm}\tau^2[e^{\tau\psi_{\eps}(\cdot,0)}u_{\pm}(\cdot,0)]^2_{1/2,\mathbb{R}^{n-1}}+\sum_{\pm}[D(e^{\tau\psi_{\eps,\pm}}u_{\pm})(\cdot,0)]^2_{1/2,\mathbb{R}^{n-1}}
\end{aligned}
\]
and 
\begin{equation*}
R_1:= (\varepsilon\tau)^{1/2}\sum_{\pm}\int_{\mathbb{R}^{n-1}}e^{2\tau\psi_{\varepsilon}(x',0)}(|D_{x_n}u_{\pm}(x',0)|^2+|D_{x'}u_{\pm}(x',0)|^2+\tau^2|u_{\pm}(x',0)|^2)dx'.
\end{equation*}
Remind that $\eta_{g,\mu}$ is defined in \eqref{etagmu}. Notice that $\Xi$ in (4.25) of \cite{dflvw} corresponds to $LHS$ here.

As in \cite[Section 4.3]{dflvw}, we introduce some local differential operators that only depend on $x_n$, in such a way that we can apply estimate \eqref{8.241}. Let us define
\begin{equation}\label{7.7}
A^\delta_{\pm}(x',x_n):=A_{\pm}(\delta x',\delta x_n),
\end{equation}
\begin{equation}\label{7.8}
\mathcal{L}_\delta(x',x_n,D)u:=\sum_{\pm}H_{\pm}{\rm div}(A^\delta_{\pm}(x',x_n)\nabla u_{\pm}),
\end{equation}
and the transmission conditions
\begin{equation*}
\begin{cases}
\theta_0(x')=u_+(x',0)-u_-(x',0),\\
\theta_1(x')=A^\delta_+(x',0)\nabla u_+(x',0)\cdot e_n-A^\delta_-(x',0)\nabla u_-(x',0)\cdot e_n.
\end{cases}
\end{equation*}
Next, recalling that $x'_g=g/\mu$ and $g\in \mathbb{Z}^{n-1}$, we define
\begin{equation*}
\begin{cases}
A^{\delta,g}_{\pm}(x_n):=A^\delta_{\pm}(x'_g, x_n)=A_{\pm}(\delta x'_g,\delta x_n),\\
\mathcal{L}_{\delta,g}(x_n,D)u:=\sum_{\pm}H_{\pm}{\rm div}(A^{\delta,g}_{\pm}(x_n)\nabla u_{\pm}).
\end{cases}
\end{equation*}
We notice that $A^{\delta,g}_{\pm}(x_n)$ satisfies assumptions \eqref{elliptic1}, \eqref{elliptic2} and also the Lipschitz condition
\begin{equation*}
|A^{\delta,g}_{\pm}(\tilde x_n)-A^{\delta,g}_{\pm}(x_n)|\leq M_0\delta|\tilde x_n-x_n|.
\end{equation*}

We now apply \eqref{8.241} to each summand and add up with respect to $g\in\mathbb{Z}^{n-1}$ to obtain that
\begin{equation}\label{8.1}
\sum_{g\in \mathbb{Z}^{n-1}}LHS(u\eta_{g,\mu})\leq C \sum_{g\in \mathbb{Z}^{n-1}}(d^{(1)}_{g,\mu}+d^{(2)}_{g,\mu}+d^{(3)}_{g,\mu}),
\end{equation}
where
\begin{equation*}
\begin{aligned}
d^{(1)}_{g,\mu}=& \int_{\mathbb{R}^n}|\mathcal{L}_{\delta,g}(x_n,D)(u\eta_{g,\mu})|^2e^{2\tau\psi_{\varepsilon}(x)}dx,\\
d^{(2)}_{g,\mu}=&\tau^{3}\int_{\mathbb{R}^{n-1}}|e^{\tau\psi_{\varepsilon}(x',0)}\theta_{0;g,\mu}(x')|^2dx'+[D_{x'}(e^{\tau\psi_{\varepsilon}}\theta_{0;g,\mu})(\cdot,0)]^2_{1/2,\mathbb{R}^{n-1}},\\
d^{(3)}_{g,\mu}=&\tau\int_{\mathbb{R}^{n-1}}|e^{\tau\psi_{\varepsilon}(x',0)}\theta_{1;g,\mu}(x')|^2dx'+[e^{\tau\psi_{\varepsilon}(\cdot,0)}\theta_{1;g,\mu}(\cdot)]^2_{1/2,\mathbb{R}^{n-1}},
\end{aligned}
\end{equation*}
where we set
\begin{equation*}
\theta_{0;g,\mu}(x'):=u_+(x',0)\eta_{g,\mu}(x')-u_-(x',0)\eta_{g,\mu}(x')=\theta_0(x')\eta_{g,\mu},
\end{equation*}
\begin{equation*}
\theta_{1;g,\mu}(x'):=A^{\delta,g}_+(0)\nabla (u_+\eta_{g,\mu})\cdot e_n-A^{\delta,g}_-(0)\nabla (u_-\eta_{g,\mu})\cdot e_n.
\end{equation*}

We now proceed as in \cite[Section 4.3, pp.201-204]{dflvw} for the estimates of the terms $d^{(j)}_{g,\mu}$, $j=1,2,3$ in  \eqref{8.1}. For the sake of clarity, we show here the estimate of the term $d^{(1)}_{g,\mu}$. By \eqref{elliptic1}, \eqref{elliptic2}, \eqref{7.4}, \eqref{6.7} and \eqref{7.7} we obtain that
\begin{equation*}
\begin{aligned}
&\;\;|\mathcal{L}_{\delta,g}(x_n,D)(u\eta_{g,\mu})|\\
&\leq |\mathcal{L}_\delta(x',x_n,D)(u\eta_{g,\mu})|+|\mathcal{L}_\delta(x',x_n,D)(u\eta_{g,\mu})-\mathcal{L}_{\delta,g}(x_n,D)(u\eta_{g,\mu})|\\
&\leq \eta_{g,\mu}|\mathcal{L}_\delta(x',x_n,D)u|+C\eta_{g,\mu}\sum_\pm|A^\delta_\pm(x',x_n)-A^\delta_\pm(x'_g,x_n)||D^2u_\pm|\\&+C\chi_{Q_{\frac{2}{\mu}}(x'_g)}\sum_\pm\left(\mu|Du_\pm|+\mu^2|u_\pm|\right)
\\
&\leq \eta_{g,\mu}|\mathcal{L}_\delta(x',x_n,D)u|+C\chi_{Q_{\frac{2}{\mu}}(x'_g)}\sum_\pm\left(\delta\mu^{-1}|D^2u_\pm|+\mu|Du_\pm|\,+\mu^2|u_\pm|\right),
\end{aligned}
\end{equation*}
which, together with \eqref{cardAg} and since $\mu=(\eps\tau)^{1/2}>1$, implies
\begin{equation}\label{8.3}
\sum_{g\in \mathbb{Z}^{n-1}}d^{(1)}_{g,\mu}\leq
C\int_{\mathbb{R}^n}|\mathcal{L}_\delta(x',x_n,D)u|^2\,e^{2\tau\psi_{\varepsilon}}dx+CR_2,
\end{equation}
where
\begin{equation*}
\begin{aligned}
R_2=\frac{\delta^2}{\mu^{2}}\sum_{\pm}\int_{\mathbb{R}^n_{\pm}}|D^2u_{\pm}|^2\,e^{2\tau\psi_{\varepsilon,\pm}}dx+\mu^{2}\sum_{\pm}\int_{\mathbb{R}^n_{\pm}}|Du_{\pm}|^2\,e^{2\tau\psi_{\varepsilon,\pm}}dx
+\mu^{4}\int_{\mathbb{R}^n}|u|^2\,e^{2\tau\psi_{\varepsilon}}dx.
\end{aligned}
\end{equation*}
With similar calculations, which are explicitly written in the above mentioned pages of \cite{dflvw}, we can estimate $d^{(2)}_{g,\mu}$, $d^{(3)}_{g,\mu}$ and get
\begin{eqnarray}\label{4.52}
LHS(u)&\leq &C\left(\int_{\mathbb{R}^n}|\mathcal{L}_\delta(x',x_n,D)u|^2\,e^{2\tau\psi_{\varepsilon}}dx+[e^{\tau\psi_{\varepsilon}(\cdot,0)}\theta_1]^2_{1/2,\mathbb{R}^{n-1}}\right.\nonumber\\
&&+[D_{x'}(e^{\tau\psi_{\varepsilon}}\theta_0)(\cdot,0)]^2_{1/2,\mathbb{R}^{n-1}}+\tau^{3}\int_{\mathbb{R}^{n-1}}e^{2\tau\psi_{\varepsilon}(x',0)}|\theta_0(x')|^2dx'\nonumber\\
&&\left.+\tau\int_{\mathbb{R}^{n-1}}e^{2\tau\psi_{\varepsilon}(x',0)}|\theta_1(x')|^2dx'+R_3\right).
\end{eqnarray}
where
\begin{eqnarray*}
R_3&=
& \frac{\delta^2}{\mu^{2}}\sum_{\pm}\int_{\mathbb{R}^n_{\pm}}|D^2u_{\pm}|^2\,e^{2\tau\psi_{\varepsilon,\pm}}dx+\mu^{2}\sum_{\pm}\int_{\mathbb{R}^n_{\pm}}|Du_{\pm}|^2\,e^{2\tau\psi_{\varepsilon,\pm}}dx\nonumber\\
&&+\mu^{4}\int_{\mathbb{R}^n}|u|^2\,e^{2\tau\psi_{\varepsilon}}dx+(\mu+\delta^2\eps^{-1})\sum_{\pm}\int_{\mathbb{R}^{n-1}}|D{u}_{\pm}(x',0)|^2e^{2\tau\psi_{\varepsilon}(x,0)}dx\nonumber\\
&&+\mu\tau^2\sum_{\pm}\int_{\mathbb{R}^{n-1}}|u_{\pm}(x',0)|^2e^{2\tau\psi_{\varepsilon}(x',0)}dx'+(\mu^4+\delta^{2}\mu^{-2}\tau^2)\sum_{\pm}[e^{\tau\psi_{\varepsilon}(\cdot,0)}u_{\pm}(\cdot,0)]^2_{1/2,\mathbb{R}^{n-1}}\nonumber\\
&&+\delta^{2}\mu^{-2}\sum_{\pm}[D(u_{\pm}e^{\tau\psi_{\varepsilon,\pm}})(\cdot,0)]^2_{1/2,\mathbb{R}^{n-1}}.
\end{eqnarray*}
We now set $\varepsilon=\delta$ and choose a sufficiently small $\delta_0$  and a sufficiently large $\tau_0$, both depending on $\lambda_0$, $\Lambda_0$, $M_0$, and $n$  such that if  $\varepsilon=\delta\leq\delta_0$ inequalities \eqref{deltaep} are satisfied and if $\tau\geq\tau_0$, then $R_3$ on the right hand side of \eqref{4.52} can be absorbed by $LHS(u)$.
We finally get the estimate \eqref{8.241} by the standard change of variable $u(\delta x',\delta x_n)$.

\section*{Acknowledgement}
Wang was partially supported by MOST 108-2115-M-002-002-MY3.

\end{document}